%% file: main.tex
\documentclass[11pt]{article}
\usepackage[utf8]{inputenc}
\usepackage[bbgreekl]{mathbbol}
\usepackage{amsmath,amssymb}
\usepackage{graphicx}
\usepackage{hyperref}
\usepackage{url}
\usepackage[numbers,square]{natbib}
\usepackage{multirow}
\usepackage{tikz}
\usepackage{tikz-cd}
\usetikzlibrary{arrows,matrix}
\usepackage{bbm}
\usepackage{authblk}
\usepackage{subcaption}

\usepackage{stmaryrd}
\usepackage{comment}
\usepackage{soul}

\usepackage[ruled,vlined]{algorithm2e}
\usepackage{float}
\usepackage{xcolor}

\usetikzlibrary{calc,patterns,shapes.geometric}
\usetikzlibrary{arrows,shapes,positioning}
\usetikzlibrary{decorations.markings,decorations.pathmorphing,decorations.pathreplacing}
\usetikzlibrary{decorations.text,decorations.shapes}

\newcommand{\la}[1]{\hat{#1}} 
\newcommand{\ca}[1]{\Bar{#1}} 

\newcommand{\R}{\mathbb{R}}

\title{A split-step Active Flux method for the Vlasov--Poisson system}

\author[1]{L. Hensel}
\author[2]{G. Gr\"unwald}
\author[1]{K. Kormann}
\author[2]{R. Grauer}

\affil[1]{Numerical Mathematics, Ruhr-Universit\"at Bochum, Universit\"atsstra{\ss}e 150, D-44801 Bochum, Germany}
\affil[2]{Theoretical Physics I, Ruhr-Universit\"at Bochum, Universit\"atsstra{\ss}e 150, D-44801 Bochum, Germany}
\date{}

\begin{document}
\maketitle

\input{content.tex}

\bigskip\bigskip

\bibliography{bibliography}

\end{document}

%% file: content.tex
		\begin{abstract}
		Active Flux is a modified Finite Volume method that evolves additional Degrees of Freedom for each cell that are located on the interface by a non-conservative method to compute high-order approximations to the numerical fluxes through the respective interface to evolve the cell-average in a conservative way. In this paper, we apply the method to the Vlasov--Poisson system describing the time evolution of the time-dependent distribution function of a collisionless plasma. In particular, we consider the evaluation of the flux integrals in higher dimensions. We propose a dimensional splitting and three types of formulations of the flux integral: a one-dimensional reconstruction of second order, a third-order reconstruction based on information along each dimension, and a third-order reconstruction based on a discrepancy formulation of the Active Flux method. Numerical results in 1D1V phase-space compare the properties of the various methods.
		\end{abstract}
		
	
	\section{Introduction}
	\label{sec:introduction}
	\noindent In this paper, we consider the Vlasov system for a collisionless plasma, which
	is a prototypical challenging system for structure-preserving bulk coupling to
	either the Poisson equation in the electrostatic limit or to the full Maxwell
	equations.
	Most space and astrophysical plasmas can be considered as collisionless. An
	estimate of the typical mean free path in space plasmas is about 1 AU
	($\approx$ distance sun $\leftrightarrow$ earth).
	Such collision-free plasmas are described very well with the kinetic Vlasov
	equation. The Vlasov equation and the field equations (Poisson and Maxwell's
	equations) are themselves linear. The bulk coupling of both sets of equations,
	which describe the interaction of particles and the electromagnetic fields,
	transforms this into a complex nonlinear system which poses a great challenge
	for numerics and a problem of highest computational expense. Therefore, there
	is a special demand for precise but also extremely effective numerical
	algorithm for the simulation of the Vlasov system.

	\noindent The mathematical model consists of the kinetic Vlasov equation
	\begin{align}\label{eqn:Vlasov}
		\partial_t f_s + \nabla_{\boldsymbol{x}} \cdot (\boldsymbol{v} f_s) + \nabla_{\boldsymbol{v}} \cdot \left( \frac{q_s}{m_s} \left( \boldsymbol{E} + \boldsymbol{v} \times \boldsymbol{B} \right) f_s \right) = 0
	\end{align}
	where $f_s = f_s(\boldsymbol{x},\boldsymbol{v},t)$ is a probability density
	function of species $s$, which depends on time $t$, macroscopic space
	$\boldsymbol{x} \in \Omega \subset\mathbb{R}^3$ and the particle velocity $\boldsymbol{v} \in \mathbb{R}^3$. Equation
	(\ref{eqn:Vlasov}) is a transport equation in
	$(\boldsymbol{x},\boldsymbol{v})$ space. The particles are coupled indirectly
	through the electromagnetic fields.
	$\boldsymbol{B}=\boldsymbol{B}(\boldsymbol{x},t) \in \mathbb{R}^3$ is the
	magnetic field and $\boldsymbol{E}=\boldsymbol{E}(\boldsymbol{x},t) \in
	\mathbb{R}^3$ is the electric field. We denote by $q_s$ and $m_s$ the charge
	and mass of particles of species $s$ (electrons: $s=e$, ions: $s=i$).
	
	\noindent The electromagnetic field satisfies Maxwell's equation, which consist
	of the dynamic evolution equations
	\begin{align}\label{eqn:Maxwell}
		\partial_t \boldsymbol{B} & = - \nabla \times \boldsymbol{E} \\
		\partial_t \boldsymbol{E} & = c^2 \left( \nabla \times \boldsymbol{B}
		- \mu_o \boldsymbol{j} \right)
	\end{align}
	and the constraints
	\begin{align}\label{eqn:constraints}
		\nabla \cdot \boldsymbol{E} & = \frac{\rho}{\epsilon_0}, \
		\nabla \cdot \boldsymbol{B}  = 0.
	\end{align}
	Here $\mu_0$ is the vacuum permeability, $\epsilon_0$ is the
	vacuum permittivity and $c$ is the speed of light. They satisfy
	the relation $1/c^2 = \mu_0 \epsilon_0$.
	Furthermore, $\rho$ is the total charge density and $\boldsymbol{j}$ is the total
	current density computed via
	\begin{align}\label{eqn:rho_j}
		\rho := \sum_s q_s n_s, \quad
		\boldsymbol{j} := \sum_s q_s n_s \boldsymbol{u}_s.
	\end{align}
	The computation of the quantities (\ref{eqn:rho_j}) requires moments
	of $f_s$, namely
	\begin{align*}
		n_s(\boldsymbol{x},t) & := \int_{\mathbb{R}^d} f_s(\boldsymbol{x},\boldsymbol{v},t) \, d\boldsymbol{v}
		& &\mbox{particle density of species $s$}\\
		\boldsymbol{u}_s (\boldsymbol{x},t) & := \frac{1}{n_s(\boldsymbol{x},t)}\int_{\mathbb{R}^d}
		\boldsymbol{v} f_s(\boldsymbol{x},\boldsymbol{v},t) \,
		d\boldsymbol{v}  & &\mbox{bulk velocity of species $s$}.
	\end{align*}
	Thus, the solution of the kinetic Vlasov equation directly influences
	the electromagnetic fields via the source term $\boldsymbol{j}$.

	\noindent In this paper, we consider the electrostatic limit, which reduces Maxwell's equations to
	\begin{align}\label{eqn:Poisson}
		\boldsymbol{B} & = 0, \\
		\Delta \phi & = -\frac{\rho}{\epsilon_0} \; , \;\; \boldsymbol{E} = - \nabla \phi \; ,
	\end{align}
	and the Vlasov equation directly influences the electric field $\boldsymbol{E}$ via the source term $\rho$.
	
	\noindent The Vlasov--Maxwell system has been used in numerous simulations of space and
	astrophysical plasmas, as well as in simulations of laboratory experiments (see
	\cite{palmroth-etal:2018,allmann_rahn-lautenbach-grauer:2021,nishikawa-etal:2021} to name a few). \\

	\noindent The numerical simulation of the Vlasov--Maxwell system is challenging for
	several reasons:
	\begin{itemize} 
		\item The Vlasov equation is \textbf{high-dimensional} (6 dimensions in
		position and velocity phase space + time), thus numerical simulations are
		computationally very expensive.
		
		\item Although the Vlasov equation has the form of a relatively simple
		transport equation, the advection velocity $\boldsymbol{v}$ might have very
		large components, which can lead to \textbf{severe time-step restrictions}.
		
		In addition, fast electron dynamics at realistic mass ratios presents an additional major challenge.

		\item The distribution function $f_s$ is \textbf{positive and conserved} and numerical schemes should conserve these properties to yield physical results. Conservation in $f_s$ can easily be obtained by using a finite volume method or a discontinuous Galerkin method. To maintain positivity a limiter will in general be necessary \cite{leveque:2002}.
		
		\item From a conservative approximation of $f_s$ we obtain a conservative approximation of the continuity equations for the particle densities $n_s$, but not necessarily for quantities associated with higher moments.
		For example, the total energy
		\begin{align*}
			E_{\text{Vlasov--Maxwell}} &:= \sum_{s} \frac{m_{s}}{2} \int\int |\mathbf{v}|^{2} f_{s}(\mathbf{x}, \mathbf{v}, t) \mathrm{d} \mathbf{v}\mathrm{d} \mathbf{x} \\
			&+ \frac{1}{2} \int \left(\varepsilon_{0}|\mathbf{E}|^{2}+\frac{1}{\mu_{0}}|\mathbf{B}|^{2}\right)\mathrm{d} \mathbf{x}
		\end{align*}
		is constant in time. This needs to be included as an additional constraint in the discretization as, e.g., in \cite{juno-hakim-etal:2018}.
		If, for example, information about the characteristics is used in the numerical procedure, then a good numerical scheme should also guarantee that by the acceleration term $\boldsymbol{v} \times  \boldsymbol{B}$ as little energy as possible is added, since artificial heating of the distribution may occur.
		\item The numerical solution  should \textbf{satisfy the constraints} (\ref{eqn:constraints}).
	\end{itemize}
	Thus, an ideal scheme for the Vlasov--Maxwell or Vlasov--Poisson system should be
	of high spatial order with low dissipation and therefore little heating and in addition of low computational cost.
	
	This challenge has been known for many years and tremendous efforts have been made to overcome these challenges.
	Numerical methods to solve the Vlasov equation range from Eulerian to
	semi-Lagrangian,
	from finite difference \cite{qin-shu:2011},
	finite volumes \cite{filbet-sonnendruecker-bertrand:2001},
	discontinuous Galerkin (DG) \cite{rossmanith-seal:2011,juno-hakim-etal:2018}
	to spectral methods \cite{eliasson:2003,delzanno:2015}.
	In addition to these grid-based methods, particle approaches such as
	particle-in-cell (PIC) \cite{arber-bennett-etal:2015}
	or fast summation techniques \cite{masek-gibbon:2010}
	are well established.
	
	All these methods have their own strengths and weaknesses that determine their
	applicability to a particular problem of plasma physics. Lagrangian schemes for
	the Vlasov equation offer the advantage that the characteristics can be easily
	calculated. 
	In particular, it is most common to split the $\textbf{x}$-advection and the $\textbf{v}$-advection parts which yields explicitly solvable characteristic equations in cases where the magnetic field can be neglected \cite{cheng-knorr:1976}.  
	Unfortunately, large splitting errors may occur in the presence of magnetic fields as demonstrated in
	\cite{schmitz-grauer:2006b} (see also \cite{roe:2017}). To circumvent these
	splitting errors, the back-substitution method was introduced
	\cite{schmitz-grauer:2006b,schmitz-grauer:2006c}. In summary, the major issues in
	Lagrangian schemes are \emph{large stencils} for high-order schemes to minimize
	heating and the application of dimensional \emph{splitting}.
	
	Eulerian schemes are generally preferable if a collision term is included. A
	comparison of Eulerian and Lagrangian schemes is presented in
	\cite{filbet-sonnendruecker:2003}. Eulerian schemes do not take full advantage
	of the simple advection form of the Vlasov equation. Thus, dissipation is
	typically slightly enhanced compared to Lagrangian schemes which leads to unphysical
	heating of the system.
	
	In this work we will apply different approaches to the Active Flux method
	\cite{eymann-roes:2011,eymann-roe:2013,roe:2021, he2020treatment}  for solving the
	Vlasov equation. We try to show the advantages, but also the necessities of
	combining the Active Flux method with splitting techniques.
	The Active Flux method has originally been introduced as a third-order finite volume method for hyperbolic conservation laws. It evolves additional point values at the cell interfaces, from which the fluxes are computed. Through this, the higher order is achieved without having to enlarge the stencil. In fact, by adding Degrees of Freedom (DOF) to the inside of the cell, Active Flux schemes of arbitrary order of accuracy  can be designed \cite{abgrall-etal:2023b}.
	Also, different strategies for achieving conservation have been proposed in the past. The original scheme \cite{eymann-roes:2011} introduced a bubble function to enforce a strict local conservation law while another approach has been discussed in \cite{he2020treatment} where the discrepancy between the cell averages from the reconstruction and those computed by enforcing the local conservation law get distributed to the nodal values. In this paper we will refer to the latter one as the discrepancy distribution formulation of Active Flux.
	A common feature is that the interface values are shared between adjacent cells, resulting in a continuous reconstruction of the solution. Active Flux is attractive in cases where the point values can be updated using exact or approximate evolution formulas, which is the case for the Vlasov system. 
	Over the last years, Active Flux has gained a lot of interest mainly
	in the fluid community \cite{barsukov-holm-etal:2019,barsukow:2021,chudzik-helzel-kerkmann:2021,bai-roe:2021}.
	In \cite{kiechle2023active}, an Active Flux method without dimensional splitting was proposed for the 1D1V Vlasov--Poisson system and a limiter introduced in \cite{kiechleNew}.
	Another important feature of Active Flux is that it is a truly
	multidimensional approach. Therefore, the most obvious idea is to apply the
	Active Flux method directly to the 6-dimensional Vlasov equation. However, this
	poses a major computational challenge: i) one needs to construct a polynomial
	reconstruction based on all the point values at the cell interfaces ($3^6 - 1 = 728$
	points in six dimensions) contrasting $2\times (3^3 -1 = 26)$ points in the case of
	splitting advection in configuration and velocity space, and ii) the
	evaluation of the six-dimensional characteristics in phase space requires
	additional intermediate calculations of the full Vlasov system to obtain the
	electromagnetic fields at the required intermediate times. Both of these
	problems are manageable as numerical methods, but would slow down any realistic
	simulation (e.g., collisionless reconnection) to the point that it is not
	feasible at all.
	Thus in this paper, we make an approximation of applying splitting techniques
	separating advection in configuration and velocity space. At that point, it
	should be noted that the advection in velocity space does not change the charge
	distribution $\rho$ and thus keeps the electric field $\boldsymbol{E}$ constant
	in time during the advection in velocity space. In this paper, we focus
	on second- and fourth-order splitting techniques. The first one is the widely
	used Strang splitting (see \cite{cheng-knorr:1976}) of the advection in
	configuration and velocity space. This simplification comes with the drawback
	of lowering the order of the scheme to only second order in time, but third
	order is maintained in space with a suitable implementation of the fluxes. If we accept the splitting of advection into
	configuration and velocity space, we are dealing with two advection problems,
	each with constant in time and space advection. To achieve an overall third-order
	accuracy, we utilize the fourth-order operator splitting described in
	\cite{yoshida1990construction}. Similar to the second-order splitting, the
	subproblems reduce to the case of a constant in time and space advection
	problem.
	We note that splitting of configuration and velocity space evolution is the common
	strategy for kinetic simulations of collisionless plasmas. It is the standard
	approach in PIC simulations (see the excellent review
	\cite{marcowith-ferrand-etal:2020} and references therein), but also used in
	spatially high-order Lagrangian Vlasov simulations \cite{kormann-reuter-rampp:2019}.
	This simplification now allows the characteristics in configuration and
	velocity space to be computed analytically. Still, the combination of Active Flux with the splitting approach is not unique, and three strategies
	will be tested in Sec.~\ref{sec:splittingAF}.
	The non-uniqueness arises from the
	necessity to switch consistently between point values, edge and face
	averages. 
	
	The paper is organized as follows: 
	In Sec.~\ref{sec:AF}, we review the original Active Flux method as well as the discrepancy distribution formulation for the linear advection problem in one dimension. 
	Sec.~\ref{sec:splittingAF} describes our splitting strategies proposed for the Vlasov--Poisson system by introducing higher-dimensional homogeneous and inhomogeneous grids that can be updated directionally by consecutive one-dimensional updates of the described approaches to the Active Flux method. Thereafter, we show how to apply directional operator splitting schemes to achieve high-order convergence in space and time.   
	The consistent solution of the Poisson equation on the introduced grid structures is discussed in Sec.~\ref{sec:poisson}.
	Sec.~\ref{sec:applications} summarizes our numerical experiments for two benchmark problem for electrostatic plasma dynamics. We compare the convergence behavior as well as the long-time behavior of the various Active flux schemes to a third-order conservative semi-Lagrangian scheme.
	Finally, Sec.~\ref{sec:conclusions} summarizes the conclusions from this comparative study.
	
	\section{Active Flux Methods}\label{sec:AF}
	\noindent 
	A numerical solver for the linear advection equation is an important component of every Vlasov solver. 
	Due to the inherently high dimensionality of the problem, the approach to increase the grid-resolution to achieve better quality solutions is generally very limited. 
	Hence, the development of new high-order methods with minimal dissipation for the Vlasov systems is necessary. 
	One relatively new type of third-order method for hyperbolic conservation laws is the here considered Active Flux method \cite{eymann-roes:2011}. 
	The general idea lies in the introduction of additional pointlike DOF that are located on the cell boundaries. 
	These point values are updated in a non-conservative fashion by a method that captures the physics of the underlying system best; for the linear advection equation this is a characteristic tracing algorithm. 
	Based on these evolved interface values numerical fluxes can be obtained by time-integration that enable for a conservative Finite Volume update of the actual cell average. \\
	The introduction of additional DOF bears several advantages. It allows for a reconstruction polynomial that is continuous among neighboring cells which makes the solution of a Riemann problem obsolete. 
	Hence, the method can be constructed to be truly multidimensional.
	Additionally, the local numerical stencil that can be achieved with the Active Flux method increases the parallelization performance on distributed memory architectures compared to finite volume methods due to the reduction of generally slow MPI communication. 
	
	\subsection{Original Active Flux Method in One Dimension} \label{sec:AF_Bubble1D}
	In this section we review the original one-dimensional Active Flux scheme as introduced by Eymann and Roe in \cite{eymann-roes:2011}. 
	We consider the hyperbolic linear one-dimensional advection equation for the
	conserved quantity $u(x,t) \in \mathbb{R}$ with the advection velocity $a \in
	\mathbb{R}$.
	\begin{align}
		\partial_t u(x,t) + \partial_x (a u(x,t)) = 0
		\label{eqn:LinearAdvection_1D}
	\end{align}
	The solution to this equation is discretized into the two types of accessible DOF namely pointlike values of the conserved quantity $u$ located on the cell boundary $u_{i+1/2}$ as well as the actual cell average of the ith cell $\hat{u}_i$
	\begin{align}
		u_{i+1/2}^n &\approx u(x_{i+1/2}, t_n) \\
		\la{u}_{i}^n &\approx \frac{1}{\Delta x} \int_{x_{i-1/2}}^{x_{i+1/2}} u(x,t^n) dx \; ,
	\end{align} 
	where we assume a uniform spacing with $\Delta x = x_{i+1/2} - x_{i-1/2}$ to simplify notations.
	Due to the fact that the cell interface points $u_{i+1/2}$ are shared among neighboring cells the number of cell internal DOF for the
	one-dimensional Active Flux method is two. The conservative update step for the
	cell average is carried out by a Finite Volume approach of the form 
	\begin{align}
		\la{u}_i^{n+1} =  \la{u}_i^n + \frac{\Delta t}{\Delta x} [F_{i-1/2} - F_{i+1/2}] \; .
	\end{align} 
	Here $F_{i+1/2}$ denotes the numerical fluxes or approximated time-averaged
	fluxes through the interface, which are given in the case of the constant
	velocity advection by
	
	\begin{align}
		F_{i+1/2} &= \frac{a}{6} (u_{i+1/2}^n + 4u_{i+1/2}^{n+1/2} + u_{i+1/2}^{n+1}) \\ 
		&\approx \frac{1}{\Delta t} \int_{t^n}^{t^{n+1}} au(x_{i+1/2},t) dt \;,
		\label{eqn:FluxIntegration_1D}
	\end{align}
	i.e.~the integrals of the fluxes are approximated by Simpson's rule. \\
	From the interface and cell average values at time $t=t^n$ the following continuous piecewise quadratic polynomial can be constructed
	\begin{align}
		q_i(\zeta) = u_{i-1/2} (\zeta-1)(2\zeta-1) + (6\la{u}_i - u_{i-1/2} - u_{i+1/2}) \zeta(1-\zeta) + u_{i+1/2} \zeta (2\zeta-1)
		\label{eqn:Reconstruction_1D}
	\end{align}     
	with $\zeta = (x-x_{i-1/2})/\Delta x$ and $0 \leq \zeta \leq 1$. The
	non-conservative update of the interface values is then carried out by tracing
	the solution back to its characteristic origin that may be located in a
	neighboring cell
	\begin{align}
		u_{i+1/2}^{n+1} = \begin{cases}
			q_{i}(\zeta_R -\frac{a \Delta t}{\Delta x}) & a \geq 0 \\
			q_{i+1}(\zeta_L -\frac{a \Delta t}{\Delta x}) & a < 0
		\end{cases}
		\label{eqn:TracingAdvection_1D} \; .
	\end{align}
	Here $\zeta_R=1$ corresponds to the right interface value and $\zeta_L=0$ to the
	left. 
	The analogous methodology can be applied to obtain $u_{i+1/2}^{n+1/2}$.
	The Courant number $\nu = (a\Delta t)/\Delta x$ is hereby restricted by the
	CFL condition $|\nu| \leq 1$ for stability. For the considered case
	of a constant velocity advection we can carry out the mentioned steps easily by
	explicit calculation, which yields the following update equations for interface
	\begin{align}
		u_{i+1/2}^{n+1} = \begin{cases}
			\nu(3\nu-2) u_{i-1/2}^n + 6\nu (1-\nu) \la{u}_i^n + (1-\nu)(1-3\nu) u_{i+1/2}^n & a \geq 0 \\
			(1-\nu)(1-3\nu)u_{i-1/2}^n + 6\nu (1-\nu) \la{u}_i^n + \nu(3\nu-2)u_{i+1/2} & a < 0
		\end{cases}
		\label{eqn:UpdateFormula_1D_Interface}
	\end{align}
	and cell average values
	\begin{align}
		\la{u}_i^{n+1} = \begin{cases}
			\nu^2 (\nu-1) u_{i-3/2}^n + \nu^2(3-2\nu) \la{u}_{i-1}^n + \nu(1-\nu) u_{i-1/2}^n + (1-\nu)^2 (1+2\nu) \la{u}_{i}^n - \nu(1-\nu)^2 u_{i+1/2} & a \geq 0 \\
			\nu^2 (\nu-1) u_{i+3/2}^n + \nu^2(3-2\nu) \la{u}_{i+1}^n + \nu(1-\nu) u_{i+1/2}^n + (1-\nu)^2 (1+2\nu) \la{u}_{i}^n - \nu(1-\nu)^2 u_{i-1/2} & a < 0 \; .
		\end{cases}
		\label{eqn:UpdateFormula_1D_Average}
	\end{align}
	\subsection{Discrepancy Distribution in One Dimension}
	Another type of Active Flux method that was introduced in
	\cite{he2020treatment} considers the cell only to be made of pointlike values
	that are updated according to the tracing algorithm and then corrected in such
	a way that global conservation is ensured. This \emph{discrepancy
		distribution} approach was initially motivated by the need for a cell-internal
	\textit{Bubble Function} for multidimensional applications to guarantee third-order
	accuracy and local conservation. 
	In the following we shortly summarize the discrepancy distribution method as it was introduced in \cite{he2020treatment}.  \\ 
	Assuming now a pointlike center value in every cell
	\begin{align}
		u_i^n = u(x_i, t^n)
	\end{align}
	one can again reconstruct via a quadratic polynomial 
	\begin{align}
		q_i(\xi) = \frac{\xi(\xi-1)}{2} u_{i-1/2} + (1-\xi^2) u_i + \frac{\xi(\xi+1)}{2} u_{i+1/2}
		\label{eqn:Reconstruction_Discrepancy_1D}
	\end{align}
	where now $-1 \leq \xi \leq 1$ determines the cell internal location. 
	The characteristic tracing method can then be
	applied the same way as shown in \eqref{eqn:TracingAdvection_1D} as a predictor
	step for $\xi_R=1$ and $\xi_L=-1$ to obtain the values $u^{*}_{i+1/2}$ and
	$u^{*}_{i-1/2}$ at the time stages $t^{n+1}$ and $t^{n+1/2}$. The center
	value at $\xi_C=0$ is also to be updated the same non-conservative way using
	the local cell reconstruction
	\begin{align}
		u(x_i, t^{n+1}) \approx u_i^{*n+1} = q_i(\xi_C - \eta), \quad u_{i+1/2}^{*n+1} = \begin{cases}
				q_{i}(\xi_R - \eta) & a \geq 0 \\
				q_{i+1}(\xi_L - \eta) & a < 0
		\end{cases}
	\end{align}
	with $\eta=(2a\Delta t)/\Delta x = 2 \nu$. With numerical fluxes obtained from time
	averaging as in eq. \eqref{eqn:FluxIntegration_1D} using now the predictor-values $u_{i+1/2}^{*n+1/2}$ and $u_{i+1/2}^{*n+1}$ we can compute the cell internal
	discrepancy $\delta u_i$ after one full time-step update as
	\begin{align}
		\delta u_i = \la{u}_i^n + \frac{\Delta t}{\Delta x} [F_{i-1/2} - F_{i+1/2}] - \la{u}_i^{n+1}
		\label{eqn:Discprepancy_1D}
	\end{align} 
	where the cell averages are approximated by
	
	\begin{align}
		\la{u}_i^n \approx \frac{u_{i-1/2}^{n} + 4u_{i}^n + u_{i+1/2}^n}{6} \\
		\la{u}_i^{n+1} \approx \frac{u_{i-1/2}^{*n+1} + 4u_{i}^{*n+1} + u_{i+1/2}^{*n+1}}{6} \; . 
	\end{align}
	One can now distribute the cell internal discrepancy in such a way that the
	corrected cell average is indeed exact. By introducing the parameters $\alpha$
	and $\beta_1, \beta_2$ corresponding to the distribution of the discrepancy to
	the center and interface values, respectively, one finds the following
	restriction that needs to be fulfilled
	\begin{align}
		\frac{(\beta_1 + \beta_2)}{6} + \frac{2 \alpha}{3} = 1. 
	\end{align} 
	Generally we symmetrically set $\beta_1 = \beta_2 = \beta$. The discrepancy is
	then distributed to the point values from the initial predictor step according
	to the parameters like
	\begin{align}
		u_i^{n+1} = u_i^{*n+1} + \alpha \delta u_i \\
		u_{i+1/2}^{n+1} = u_{i+1/2}^{*n+1} + \beta \frac{\delta u_i + \delta u_{i+1}}{2}.
	\end{align} 
	Choosing the parameters as $\alpha=\beta=1$ means the equal distribution of the
	discrepancy on all accessible DOF in the cell. By this method the clear distinction between
	cell averages and interface values that was present in the previous section (\ref{sec:AF_Bubble1D}) can
	be dissolved. Hence, after distributing the discrepancy over the whole cell all
	accessible DOF posses the same interpretation after correction. 
	The CFL condition is furthermore more restrictive than before as we need to have $|\nu| = |a \Delta t / \Delta x| \leq 1 / 2$ to ensure stability for the tracing update of the center values.  \\
	We can again derive analytical update formulae for the discrepancy distribution Active Flux scheme.
	For positive advection speeds $a \geq 0$ the update of the individual cell DOF can be carried out as
	
	\begin{align}
		u_{i+1/2}^{n+1} &= \frac{\beta \eta}{48} (\eta - 2)(2\eta -1) u_{i-3/2}^n \label{eqn:DiscrepancyFormulaInterfacePosVel}
		- \frac{\beta \eta}{12} (\eta^2 - 4\eta +2) u_{i-2}^n \\
		&+ \frac{\eta}{48} (2 \beta \eta^2 - 23 \beta \eta + 14 \beta + 24\eta - 24) u_{i-1/2}^n \notag \\ 
		&+ \frac{\eta}{6} (3\beta \eta - 2\beta - 6\eta + 12) u_{i}^n 
		- \frac{1}{48} (2\beta \eta^3 + 19\beta \eta^2 - 14\beta \eta - 24 \eta^2 + 72\eta -48) u_{i+1/2}^n \notag \\
		&+ \frac{\beta \eta}{12} (\eta^2 + 2\eta -2) u_{i+1}^n
		- \frac{\beta \eta}{48} (2\eta^2 + \eta -2) u_{i+3/2}^n  \notag \\
		u_i^{n+1} &=  \frac{1}{6} (\alpha \eta^3 + 2\alpha \eta^2 - 2\alpha \eta - 6\eta^2 + 6) u_i^n 
		- \frac{1}{24}(2 \alpha \eta^3 + \alpha \eta^2 - 2\alpha \eta - 12\eta^2 + 12\alpha) u_{i+1/2}^n \\
		&- \frac{1}{4}(3 \alpha \eta^2 - 2\alpha \eta - 2\eta^2 - 2\eta) u_{i-1/2}^n 
		- \frac{1}{6}(\alpha \eta^3 - 4 \alpha \eta^2 + 2\alpha \eta) u_{i-1}^n \notag \\
		&+ \frac{1}{24}(2\alpha \eta^3 - 5 \alpha \eta^2 + 2\alpha \eta) u_{i-3/2}^n	\notag	
	\end{align}
	and similarly for $a < 0$:
	\begin{align}
		u_{i+1/2}^{n+1} &= \frac{\beta \eta}{48} (\eta - 2)(2\eta -1) u_{i+3/2}^n 
		- \frac{\beta \eta}{12} (\eta^2 - 4\eta +2) u_{i+2}^n \\
		&+ \frac{\eta}{48} (2 \beta \eta^2 - 23 \beta \eta + 14 \beta + 24\eta - 24) u_{i+1/2}^n \notag \\
		&+ \frac{\eta}{6} (3\beta \eta - 2\beta - 6\eta + 12) u_{i}^n  
		- \frac{1}{48} (2\beta \eta^3 + 19\beta \eta^2 - 14\beta \eta - 24 \eta^2 + 72\eta -48) u_{i-1/2}^n \notag \\
		&+ \frac{\beta \eta}{12} (\eta^2 + 2\eta -2) u_{i-1}^n
		- \frac{\beta \eta}{48} (2\eta^2 + \eta -2) u_{i-3/2}^n \notag \\
		u_i^{n+1} &=  \frac{1}{6} (\alpha \eta^3 + 2\alpha \eta^2 - 2\alpha \eta - 6\eta^2 + 6) u_i^n \label{eqn:DiscrepancyFormulaCenterNegVel}
		- \frac{1}{24}(2 \alpha \eta^3 + \alpha \eta^2 - 2\alpha \eta - 12\eta^2 + 12\alpha) u_{i-1/2}^n \\
		&- \frac{1}{4}(3 \alpha \eta^2 - 2\alpha \eta - 2\eta^2 - 2\eta) u_{i+1/2}^n 
		- \frac{1}{6}(\alpha \eta^3 - 4 \alpha \eta^2 + 2\alpha \eta) u_{i+1}^n \notag \\
		&+ \frac{1}{24}(2\alpha \eta^3 - 5 \alpha \eta^2 + 2\alpha \eta) u_{i+3/2}^n	\notag	
	\end{align}
	
	\section{Splitting Strategies for the Active Flux Method in Higher Dimensions}\label{sec:splittingAF}
	Originally, an important motivation for the introduction of the Active Flux method was to develop a truly multidimensional approach for hyperbolic systems. This has been successfully implemented for various two- and three-dimensional problems such as advection, Euler equations, shallow water and acoustic equation, to name a few
	\cite{chudzik-etal:2024, abgrall-etal:2023}.
	However, a drawback of the additional interface DOF in the Active Flux method is the inefficient scaling of memory cost when extending to higher-dimensional systems. This is evident when considering the straightforward extension of the original Active Flux method, as discussed in section \ref{sec:AF_Bubble1D}, to a multidimensional Cartesian grid approach, as presented in \cite{maeng2017advective}.
	It can be seen that the number of total cell DOF scales as $n_{tot} \propto 3^d$ neglecting the domain boundaries, where $d$ is the dimensionality of the grid cell.
	Since most of these values are on the cell interfaces, they are effectively shared between neighboring cells.
	The effective number of DOF that need to be stored in memory for each cell individually reduces thereby for the cartesian grid to $n_{ind} \propto 2^d$.         
	However since we are aiming to develop strategies for applying the Active Flux method to the full six-dimensional (three spatial / three velocity directions) Vlasov system, an unsplit treatment is still not feasible.
	E.g. each cell in six-dimensional phase space would consist of $3^6 = 729$, for which correct tracing updates and flux integrals of equal cost and dimensionality need to be computed respectively.
	For this work  we rely on the well-established strategy of directional splitting for high-dimensional problems. The aim thereby is to apply one-dimensional operations along each direction consecutively. This reduces the computational expense to a manageable level in the sense that e.g. not all the $729$ DOF get updated in a single step but one-dimensional slices of three DOF per cell, as in Sec.~\ref{sec:AF_Bubble1D}, get updated along each direction.
	This comes at the price of an additional time-step dependent numerical error introduced by the splitting. A more detailed discussion of different time splitting schemes can be found in the upcoming Sec.~\ref{sec:splitting}. \\
	We would like to emphasize that a splitting approach is preferable for the Vlasov system Eq.~(\ref{eqn:Vlasov}) in particular, since the splitting in $x$ space is indeed exact without introducing any additional error, since the evolution operators for the $x,y,z$ advection all commute.
	Next, we discuss different grid structures that are suitable for extending one-dimensional Active Flux methods to multidimensional grids. Thereafter, we examine suitable operator splitting schemes to achieve higher order accuracy for the time integration of the Vlasov system. \\
	In the following we are going to propose several extensions of the one-dimensional Active Flux schemes as an advection solver to higher dimensional systems. For readability and the sake of simplicity we often consider a two-dimensional system as an example while extendibility to even higher dimensions is straightforward in most cases as will be denoted. 
	
	\subsection{Possible Grid Arrangements For Higher Dimensional Active Flux method} 
\begin{figure}[ht]
		\centerline{\includegraphics[scale=0.45]{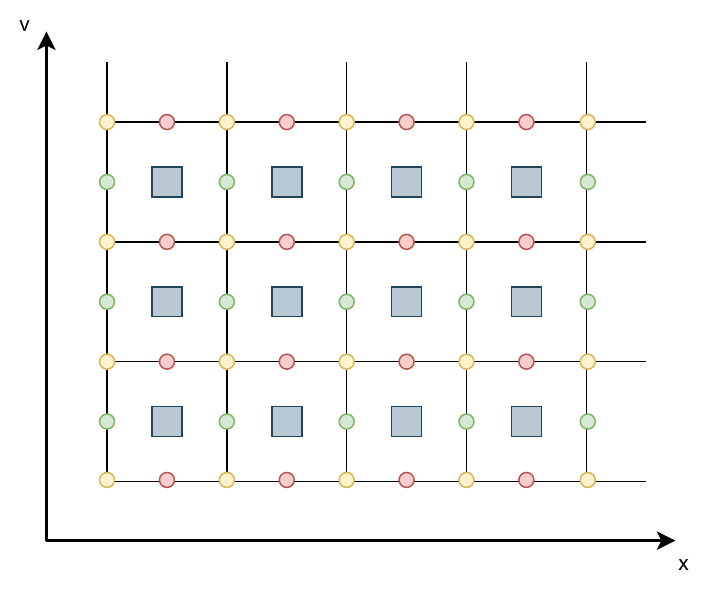}
			\includegraphics[scale=0.45]{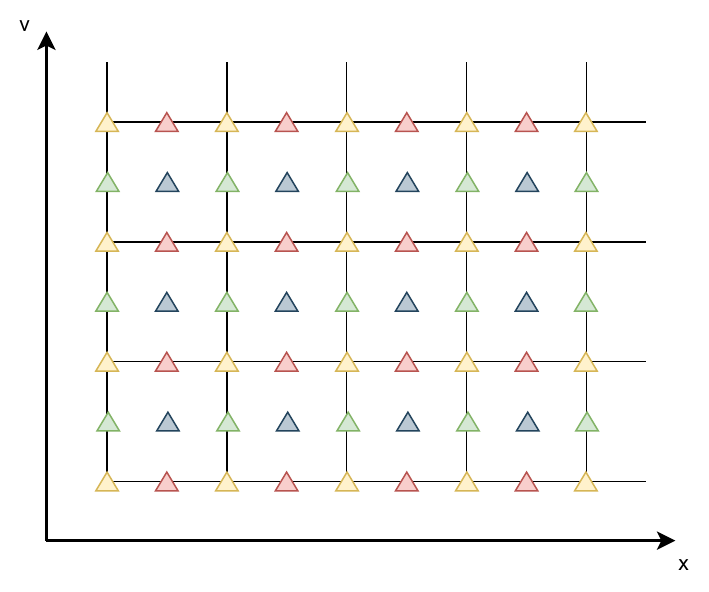}
			\includegraphics[scale=0.45]{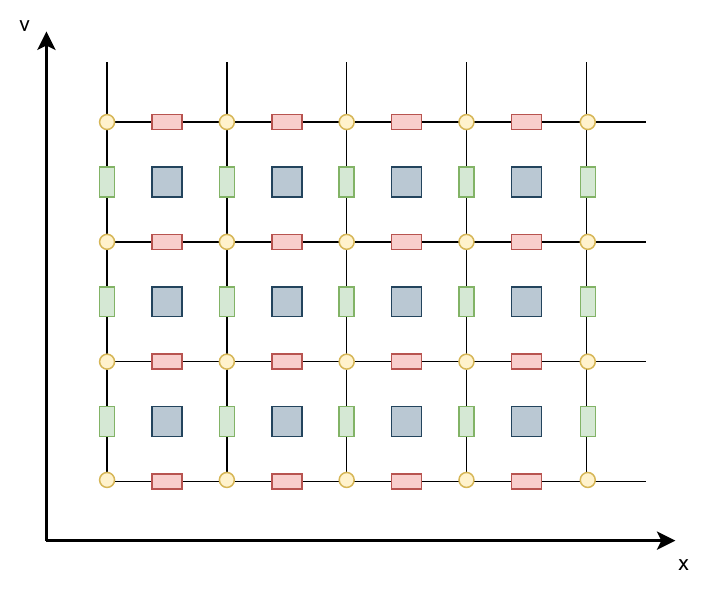}
		}
		\caption{Two-dimensional Active Flux grid arrangements:
			a) unsplit (left),  
			b) split discrepancy distribution (middle),
			c) split (right)}
		\label{fig:AF_arrangement}
	\end{figure}

	Our intention is to apply the previously described one-dimensional procedures in the context of kinetic plasma simulations to update the two-dimensional distribution function in $xv$-phase space. 
	In Fig.~\ref{fig:AF_arrangement} possible grid-arrangements for higher dimensional Active Flux method are displayed: 
	a) shows the grid for an unsplit Active Flux method in which each grid-cell is simultaneously evolved in both directions every time step. This grid consisting of cells made up of eight shared interface-point-values and a distinct cell-average is equivalent to the grid originally presented in \cite{maeng2017advective} and was already applied for the unsplit numerical treatment of the Vlasov--Poisson system in \cite{kiechle2023active}. 
	
	Since we are aiming at a directional splitting together with the described one-dimensional Active Flux advection updates we need to alter the two-dimensional grid in a way to contain the one-dimensional Active Flux cells along every direction: in b) a homogeneous computational grid suitable for the discrepancy distribution method is shown. It consists solely of point-values that are in every update step along both directions corrected in a way to give the correct cell average after integration over the respective two-dimensional grid cell. 
	In c) a grid that corresponds to the original formulation of the Active Flux method is shown. 
	The grid in this case is inhomogeneously made up from nodal point values as well as line-averages on the edges, and further integrals over all $j$-faces up to the facets and finally cell averages for each cell.
 
	As a consequence of the inhomogeneous grid we need to take into account the different types of quantities; especially when solving the Poisson equation coupled to the evolution of the distribution function in the case of Vlasov systems. \\
	As we aim to carry out exemplary simulations on the two-dimensional Vlasov system we explain the general grid setup and introduce the notation we use to distinguish the different types of quantities for the inhomogeneous grids:  
	For a two-dimensional system the grid is initially set up  by evaluating the given analytical initial condition $f_0(x, v)$ at all points $(x_i, v_j)$ with $i \in \{0, 1/2, 1, ..., N_x \}$ and $j \in \{0, 1/2, 1, ..., N_v \}$. This grid of point values is then used to initially compute the one-dimensional line- and two-dimensional cell-averages if needed e.g. by using Simpsons rule, that is easily applicable on the structure of the Active Flux grid:
	\begin{align}
		\frac{1}{\Delta v} \int_{v_{j-1/2}}^{v_{j+1/2}} f(x_{i+1/2},v,t^n) dv &\approx \frac{1}{6} (f_{i+1/2, j-1/2} + f_{i+1/2, j+1/2} + 4 f_{i+1/2, j}) =: \Hat{f}_{i+1/2, j}^n \label{eq:1D_average_v}\; , \\
		\frac{1}{\Delta x} \int_{x_{i-1/2}}^{x_{i+1/2}} f(x,v_{j+1/2},t^n) dx &\approx \frac{1}{6} (f_{i-1/2, j+1/2} + f_{i+1/2, j+1/2} + 4 f_{i, j+1/2}) =: \Hat{f}_{i, j+1/2}^n \label{eq:1D_average_x} \; , \\
		\frac{1}{\Delta x \Delta v} \int_{v_{j-1/2}}^{v_{j+1/2}} \int_{x_{i-1/2}}^{x_{i+1/2}} f(x,v,t^n)dx dv  
		&\approx \frac{1}{36} [(f_{i-1/2, j-1/2} + f_{i+1/2, j-1/2} + f_{i-1/2, j+1/2} + f_{i+1/2, j+1/2}) \notag \\
		&+ 4(f_{i, j-1/2} + f_{i, j+1/2} + f_{i-1/2, j} + f_{i+1/2, j}) + 16 f_{i, j}] =: \ca{f}^n_{i,j} \; \label{eqn:SimpsonAverage2D}
	\end{align}
	where we denote the line-averaged quantities in the grid by a hat and the cell-averaged quantities by a bar. This notation will be used consistently throughout the paper. 
	Point values can be obtained at any time by evaluating the corresponding reconstruction polynomial at a given position.   
	
	\subsection{Time Splitting Strategies}\label{sec:splitting}	

	In this subsection, we discuss the directional splitting for the Vlasov--Poisson system. The Vlasov--Poisson system can be viewed as an advection in configuration and an advection in velocity space. The advection speed for transport in velocity space is given by the electric field $E$ which is calculated from first moments of the distribution functions. This motivates a very natural splitting for the evolution of the distribution functions by separating advection in configuration and velocity space. In the six-dimensional case, the operators $L_X$ and $L_V$ can be further split into three one-dimensional advection steps. Only one Poisson solve after the three $x$-steps is necessary also in this case. 
	To be more precise, the necessary operations are introduced in Alg.~\ref{alg:Operators}.

	\SetKwComment{Comment}{/* }{ */}
	\begin{algorithm}[h!]
		\caption{Definition of operators for advection in configuration and velocity space}
		\label{alg:Operators}
		\DontPrintSemicolon
		\SetKwFunction{OStepX}{$L_X$}
		\SetKwFunction{OStepV}{$L_V$}
		\SetKwProg{Fn}{Operator}{:}{}
		\Fn{\OStepX{$f$,$\Delta t$}}{
			propagate $f$ over $\Delta t$ / solve: $\partial_t f+ \nabla_x \cdot (\boldsymbol{v} f) = 0$\; 
			calculate density: $\rho(\boldsymbol{x},t^{n+1}) = \int d^dv f(\boldsymbol{x},\boldsymbol{v},t^{n+1})$ \;
			solve Poisson equation: $\Delta \phi(\boldsymbol{x},t) = -\rho(\boldsymbol{x},\boldsymbol{v},t^{n+1})/\epsilon_0 \; , \;\; \boldsymbol{E}(\boldsymbol{x},t) = -\nabla \phi$ \; 
			return f, $\boldsymbol{E}$\;
		}
		
		\Fn{\OStepV{$f$, $\boldsymbol{E}$, $\Delta t$}}{
			propagate $f$ over $\Delta t$ / solve: $\partial_t f +  \frac{q_e}{m_e} \nabla_v \cdot  ( \boldsymbol{E}f) = 0$\;
			return f\;
		}
	\end{algorithm}

	Note that during the advection step in velocity space, indicated by $L_V$, the charge density $\rho$ and thus the electric field $E$ remain constant. Thus, both propagation steps in configuration and velocity space correspond to linear advection with constant advection in time and space.

	The accuracy of the time-integration $f^n \rightarrow f^{n+1}$ is in fact determined by the employed operator splitting scheme. The classical first order accurate splitting, called Lie splitting, is shown in Alg.~\ref{alg:Lie}.

	\begin{algorithm}[h!]
		\caption{First-order Lie splitting}
		\label{alg:Lie}
		\DontPrintSemicolon
		\While{$t \leq t_{max}$}
		{
			$f_X(\boldsymbol{x},\boldsymbol{v}, t^{n+1}), \boldsymbol{E} = L_X( f(\boldsymbol{x},\boldsymbol{v}, t^n), \Delta t)$ \; 
			$f_V(\boldsymbol{x},\boldsymbol{v}, t^{n+1}) = L_V(f_X(\boldsymbol{x},\boldsymbol{v}, t^{n+1}),\boldsymbol{E}, \Delta t) $ \; 
			solution after one timestep: $f(\boldsymbol{x},\boldsymbol{v}, t^{n+1}) = f_V(\boldsymbol{x},\boldsymbol{v}, t^{n+1})$ \;
			$t = t  + \Delta t$\;
		}
	\end{algorithm}

	Higher order schemes can be developed by carrying out additional intermediate propagation steps. 
	Using the shorthand operator notation, the second-order Strang splitting for the Vlasov--Poisson system, which consists of two half-time steps in the $x$-direction and one full-time update in the $v$-direction, is formulated in Alg.~\ref{alg:Strang}.
	\begin{algorithm}[h!]
		\caption{Second-order Strang splitting}
		\label{alg:Strang}
		\DontPrintSemicolon
		
		\While{$t \leq t_{max}$}
		{
			$f_X(\boldsymbol{x},\boldsymbol{v}, t^{n+1/2}), \boldsymbol{E} = L_X( f(\boldsymbol{x},\boldsymbol{v}, t^n), \Delta t/2)$ \; 
			$f_V(\boldsymbol{x},\boldsymbol{v}, t^{n+1}) = L_V(f_X(\boldsymbol{x},\boldsymbol{v}, t^{n+1/2}), \boldsymbol{E}, \Delta t)$ \; 
			$f_X(\boldsymbol{x},\boldsymbol{v}, t^{n+1}), \boldsymbol{E} = L_X( f_V(\boldsymbol{x},\boldsymbol{v}, t^{n+1}), \Delta t/2)$ \;
			\;
			$f(\boldsymbol{x},\boldsymbol{v}, t^{n+1}) = f_X(\boldsymbol{x},\boldsymbol{v}, t^{n+1})$ \;
			$t = t + \Delta t$
		}
	\end{algorithm}
	The time-splitting in Alg.~\ref{alg:Strang} is in fact a standard approach that was first applied to the Vlasov--Poisson system in \cite{cheng-knorr:1976}. The coupled electrostatic Poisson equation is solved before every $v$-step. 

	Higher order schemes can be constructed on the basis of these two fundamental splittings. Alg.~\ref{alg:Yoshida} shows a fourth order splitting scheme introduced in \cite{yoshida1990construction}. This method consists of three consecutive Strang splittings with specified step sizes given in terms of the constants $\gamma_1 = 2 / (2-2^{1/3})$ and $\gamma_2 = -2^{1/3} / (2-2^{1/3})$.

	\begin{algorithm}[h!]
		\caption{Fourth order Yoshida splitting}
		\label{alg:Yoshida}
		\DontPrintSemicolon
		$\gamma_1 = \frac{1}{2 - 2^{1/3}}$,  $\gamma_2 = -\frac{2^{1/3}}{2 - 2^{1/3}}$ \;	
		\While{$t \leq t_{max}$}
		{
			$f_X(\boldsymbol{x},\boldsymbol{v}, t^{n+\gamma_1/2}), \boldsymbol{E} = L_X( f(\boldsymbol{x},\boldsymbol{v}, t^n), (\gamma_1 \Delta t)/2)$ \;
			$f_V(\boldsymbol{x},\boldsymbol{v}, t^{n+\gamma_1}) = L_V(f_X(\boldsymbol{x},\boldsymbol{v}, t^{n+\gamma_1/2}),\boldsymbol{E},\gamma_1 \Delta t)$ \;
			$f_X(\boldsymbol{x},\boldsymbol{v}, t^{n+\gamma_1}), \boldsymbol{E} = L_X(( f_V(\boldsymbol{x},\boldsymbol{v}, t^{n+\gamma_1}), (\gamma_1 \Delta t)/2)$ \;
			$f(\boldsymbol{x},\boldsymbol{v}, t^{n+\gamma_1}) = f_X(\boldsymbol{x},\boldsymbol{v}, t^{n+\gamma_1})$ \;
			\;
			$f_X(\boldsymbol{x},\boldsymbol{v}, t^{n+\gamma_1 + \gamma_2 / 2}), \boldsymbol{E} = L_X( f(\boldsymbol{x},\boldsymbol{v}, t^{n + \gamma_1}), (\gamma_2 \Delta t)/2) $ \;
			$f_V(\boldsymbol{x},\boldsymbol{v}, t^{n+\gamma_1 + \gamma_2}) = L_V( f_X(\boldsymbol{x},\boldsymbol{v}, t^{n+\gamma_1 + \gamma_2/2}),\boldsymbol{E},\gamma_2 \Delta t)$ \;
			$f_X(\boldsymbol{x},\boldsymbol{v}, t^{n+\gamma_1 + \gamma_2}), \boldsymbol{E} = L_X( f_V(\boldsymbol{x},\boldsymbol{v}, t^{n+\gamma_1}), (\gamma_2 \Delta t)/2)$ \;
			$f(\boldsymbol{x},\boldsymbol{v}, t^{n+\gamma_1 + \gamma_2}) = f_X(\boldsymbol{x},\boldsymbol{v}, t^{n+\gamma_1 + \gamma_2})$ \;
			\;
			$f_X(\boldsymbol{x},\boldsymbol{v}, t^{n+\gamma_1 + \gamma_2 + \gamma_1 / 2}), \boldsymbol{E} = L_X(f(\boldsymbol{x},\boldsymbol{v}, t^{n + \gamma_1 + \gamma_2}), (\gamma_1 \Delta t)/2)$ \;
			$f_V(\boldsymbol{x},\boldsymbol{v}, t^{n+\gamma_1 + \gamma_2 + \gamma_1}) = L_V( f_X(\boldsymbol{x},\boldsymbol{v}, t^{n+\gamma_1 + \gamma_2 + \gamma_1 / 2}),\boldsymbol{E},\gamma_1 \Delta t)$ \;
			$f_X(\boldsymbol{x},\boldsymbol{v}, t^{n+\gamma_1 + \gamma_2 + \gamma_1}), \boldsymbol{E} = L_X(f_V(\boldsymbol{x},\boldsymbol{v}, t^{n+\gamma_1 + \gamma_2 + \gamma_1}), (\gamma_1 \Delta t)/2)$ \;
			\;
			$f(\boldsymbol{x},\boldsymbol{v}, t^{n+1}) = f(\boldsymbol{x},\boldsymbol{v}, t^{n+\gamma_1 + \gamma_2 + \gamma_1}) = f_X(\boldsymbol{x},\boldsymbol{v}, t^{n+\gamma_1 + \gamma_2 + \gamma_1})$ \;
			$t = t + \Delta t $ \;
		}
	\end{algorithm}
	For every splitting the Poisson equation is solved after the propagation in the physical space to determine the electrostatic potential in accordance with the updated density. 
	Our goal in this paper is to describe high-order strategies for solving the advective part of the Vlasov--Poisson system with the Active Flux method in combination with the previously described operator splitting approach for the time-integration. This means applying consecutive one-dimensional Active Flux advection updates as described in Sec.~\ref{sec:AF} to update the two-dimensional grids presented in Fig.~\ref{fig:AF_arrangement}. 
	In the upcoming subsections we are going to propose three variants for effectively using time-splitting together with the Active Flux methods exemplary for the 1D1V Vlasov problem.  
	
	\subsection{Second-Order Accurate Flux Integral}
\begin{figure}[ht]
		\centerline{\includegraphics[scale=0.65]{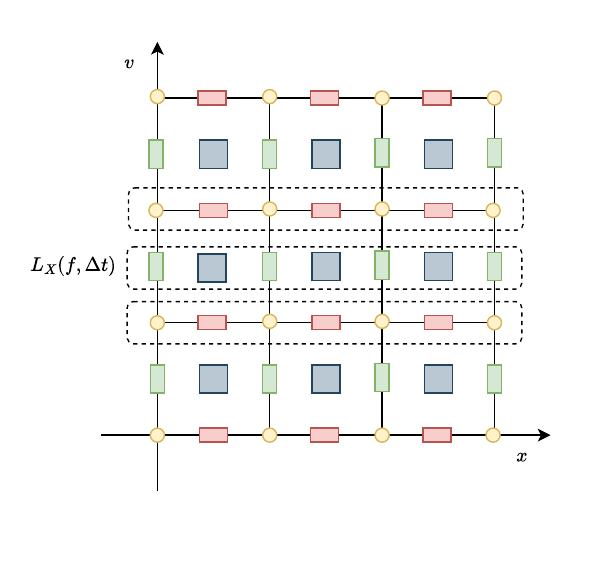}
			\includegraphics[scale=0.65]{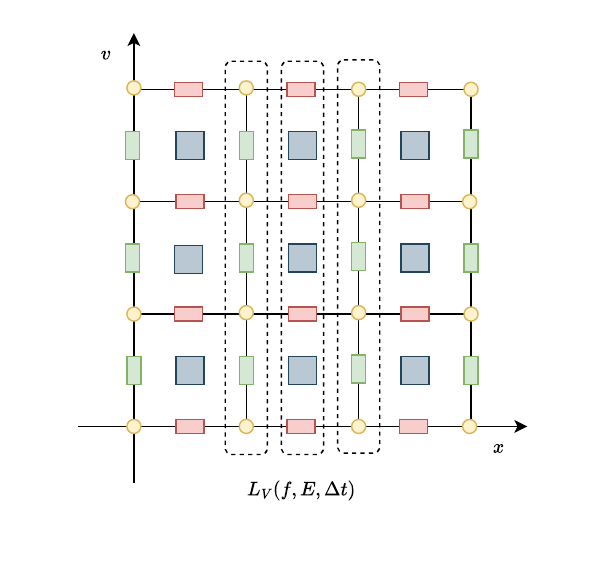}
		}
		\caption{Naive Splitting for updating every DOF of a two-dimensional Active Flux cell in phase space. This approach results in a second-order accurate approximation of the flux-integral. The evolution operators $L_X$ (left) and $L_V$ (right) are applied on one-dimensional slices of the grid along each direction. Consequently, the resulting one-dimensional Active Flux cells within these slices consist of point values on the nodes and one-dimensional line averages on the edges or line averages and two-dimensional cell averages respectively. The cell averages of any cell are conservatively updated along both directions.}
		\label{fig:2nd_order_naive}
	\end{figure}

	One straightforward strategy is the naive application of the directional one-dimensional update procedures to the original higher dimensional Active Flux grid as depicted in Fig. \ref{fig:2nd_order_naive}. \\ 
	The described naive splitting leads to an overall second-order accurate approximation of the flux integral through a given boundary and consequently a second-order update of the cell average; in this case the fluxes for the conservation update of the two-dimensional cell average are computed as time integrals over the products of previously updated line averages over the respective edge instead of actually computing the double-integral      
	\begin{align}
		&\frac{1}{\Delta v \Delta t} \int_{t^n}^{t^{n+1}} \int_{v_{i-1/2}}^{v_{i+1/2}} v \cdot f(x,v,t) dv dt \notag \\
		 &\approx \frac{1}{\Delta v \Delta t}\int_{t^n}^{t^{n+1}} \left(\int_{v_{j-1/2}}^{v_{j+1/2}} v dv \cdot \int_{v_{j-1/2}}^{v_{j+1/2}} f(x,v,t) dv \right) dt \; , \\  
		&\frac{1}{\Delta x \Delta t }\int_{t^n}^{t^{n+1}} \int_{x_{i-1/2}}^{x_{i+1/2}} E(x,t) \cdot f(x,v,t)  dx dt \notag \\
		&\approx \frac{1}{\Delta x \Delta t}\int_{t^n}^{t^{n+1}} \left(\int_{x_{i-1/2}}^{x_{i+1/2}} E(x,t) dx \cdot \int_{x_{i-1/2}}^{x_{i+1/2}} f(x,v,t) dx \right) dt \; ,
	\end{align}
	with the one-dimensional line averages (\ref{eq:1D_average_v}) and (\ref{eq:1D_average_x}) that are effectively obtained as results from one-dimensional Active Flux steps appearing on the right side of the approximation as well as the quantities
	\begin{align}	
		\frac{1}{\Delta v} \int_{v_{j-1/2}}^{v_{j+1/2}} v dv = v_j, \quad
		\frac{1}{\Delta x} \int_{x_{i-1/2}}^{x_{i+1/2}} E(x,t^n) dx \approx \Hat{E}^n_i \; .
	\end{align} 
	Find a algorithmic description of the entire update procedure in Alg.~\ref{alg:SecondOrderFluxIntegral}. 
	This version is particularly easy to implement since it only consists of one-dimensional slice-wise operations which can be straightforwardly applied to higher dimensional Vlasov systems as well. \\
	
	\begin{minipage}{0.88\linewidth}
	\begin{algorithm}[H]
		\caption{Second-Order Flux Integral Advection Operators}
		\label{alg:SecondOrderFluxIntegral}
		\DontPrintSemicolon
		\SetKwFunction{OStepX}{$L_X$}
		\SetKwFunction{OStepV}{$L_V$}
		\SetKwProg{Fn}{Operator}{:}{}
		\Fn{\OStepX{$f$,$\Delta t$}}{
			propagate $f$ over $\Delta t$ / solve: $\partial_t f+ v\partial_x f = 0$\; 
			\begin{itemize}
				\item Apply the scheme from Sec.~\ref{sec:AF_Bubble1D} / use the analytic update formulae (\ref{eqn:UpdateFormula_1D_Interface}) and (\ref{eqn:UpdateFormula_1D_Average}) to update one-dimensional horizontal slices of the grid.
				\item For slices at $v =  v_{j+1/2}$ the nodal point values $f_{i+1/2,j+1/2}$ act as interfaces and the horizontal line averages in $x$-direction $\hat{f}_{i,j+1/2}$ act as cell-averages during the update procedure.  
				\item For slices at $v =  v_{j}$ the vertical line-averages $\hat{f}_{i+1/2,j}$ act as interfaces and the two-dimensional cell-averages $\Bar{f}_{i,j}$ act as one-dimensional cell-averages during the update procedure along that direction.    
			\end{itemize}
			calculate density: $\rho(x,t^{n+1}) = \int  f(x,v,t^{n+1}) dv$ \;
			solve Poisson equation: $\partial_x^2 \phi(x,t) = -\rho(x,v,t^{n+1})/\epsilon_0 \; , \;\; E(x,t) = -\partial_x \phi$ \; 
			\begin{itemize}
				\item Apply the scheme from Sec.~\ref{sec:poisson} to compute the electric field from the distribution function on the inhomogeneous grid. 
				\item The electric field is consequently in the format: point values on the interfaces $E_{i+1/2}$ and $x$-averages $\hat{E}_i$ at the center locations.
			\end{itemize}
			return f, $E$\;
		}
		\Fn{\OStepV{$f$, $E$, $\Delta t$}}{
			propagate $f$ over $\Delta t$ / solve: $\partial_t f + (q_e / m_e) E \partial_v f = 0$\;
			\begin{itemize}
				\item Apply the scheme from Sec.~\ref{sec:AF_Bubble1D} / use the analytic update formulae (\ref{eqn:UpdateFormula_1D_Interface}) and (\ref{eqn:UpdateFormula_1D_Average}) to update one-dimensional horizontal slices of the grid.
				\item For slices at $x =  x_{i+1/2}$ the nodal point values $f_{i+1/2,j+1/2}$ act as interfaces and the horizontal line averages in $v$-direction $\hat{f}_{i+1/2,j}$ act as cell-averages during the update procedure. 
				\item For slices at $x =  x_{i}$ the horizontal line-averages $\hat{f}_{i,j+1/2}$ act as interfaces and the two-dimensional cell-averages $\Bar{f}_{i,j}$ act as one-dimensional cell-averages during the update procedure along that direction.   
			\end{itemize}
			return f\;
		}
		\;
	\end{algorithm}
	\end{minipage}
	
	\subsection{Third-Order Accurate Flux Integral}
\begin{figure}[ht]
		\centerline{\includegraphics[scale=0.55]{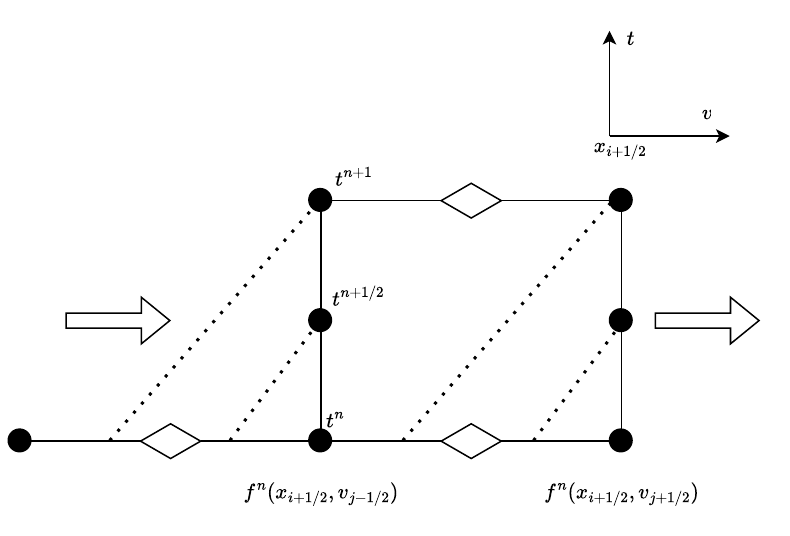}
			\includegraphics[scale=0.55]{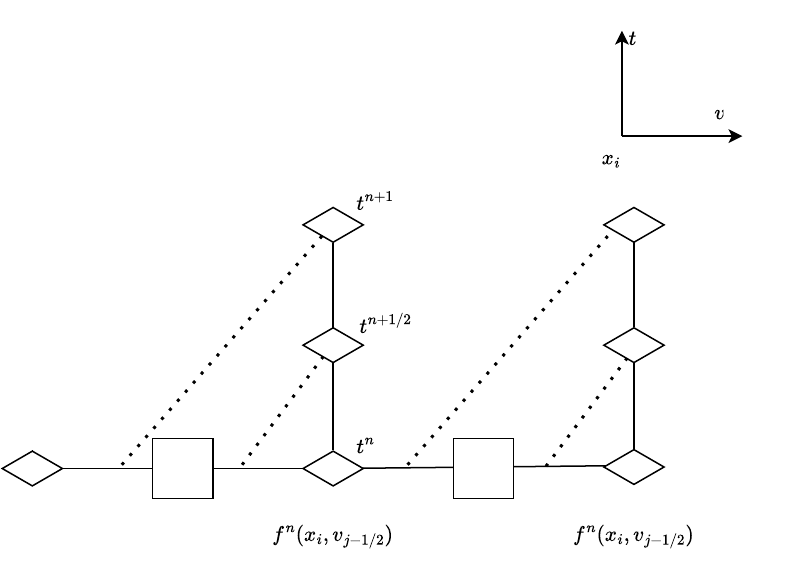}
		}
		\caption{Exemplary one-dimensional updates in the $v$-dimension for computing the third-order fluxintegral over an edge of a two-dimensional AF cell. Characteristic tracing is indicated by the dotted lines: (left) on the edges along the $v$-direction, with conservation update of the line average of that edge, (right) on the centers of the $x$-edges without conservation update of the cell average which is carried out by computing the flux integral over the entire $x$-edge in space and time.}
		\label{fig:1D_steps}
	\end{figure}

	For simplicity of notation let us first discuss the reduced case with $x \in \Omega \subset \R$ and $v \in \R$. 
	We can rewrite the Vlasov equation in terms of fluxes of the form $g(f(x,v,t)) = vf(x,v,t)$ and $h(f(x,v,t)) = E(x,t) f(x,v,t)$
	\begin{align}
		\partial_t f(x,v,t) + \partial_x g(f(x,v,t)) + \partial_v h(f(x,v,t)) = 0
	\end{align}
	Reformulating this problem by applying the Finite Volume method, namely discretizing the distribution function into averages over 2D volumes or rectangular cells $F_{i,j}^{n}$ yields the following update formula in discretized time:
	\begin{align}
		\ca{f}_{i,j}^{n+1} = \ca{f}_{i,j}^{n} - \frac{\Delta t}{\Delta x} (G_{i + 1/2, j} - G_{i - 1/2, j}) - \frac{\Delta t}{\Delta v} (H_{i, j+1/2} - H_{i, j-1/2})
	\end{align} 
	with the numerical fluxes being given by time integration over the period $[t^n, t^{n+1}]$ and over vertical- 
	\begin{align}
		G_{i+1/2, j} = \frac{1}{\Delta t \Delta v}\int_{t^n}^{t^{n+1}} \int_{v_{j-1/2}}^{v_{j+1/2}} g(f(x_{i+1/2},v,t)) dv dt
	\end{align}
	and horizontal edges in phase space respectively
	\begin{align}
		H_{i, j+1/2} = \frac{1}{\Delta t \Delta x} \int_{t^n}^{t^{n+1}} \int_{x_{i-1/2}}^{x_{i+1/2}} h(f(x, v_{j+1/2}, t)) dx dt
	\end{align}
	Through directional one-dimensional Active Flux updates on the grid Fig.~\ref{fig:AF_arrangement}.c it is possible to obtain third-order accurate cell information for the nodal point values and the one-dimensional averages at the time stages $t^n$, $t^{n+1/2}$ and $t^{n+1}$ (see Fig.~\ref{fig:1D_steps}).
	Consequently, we intend to update the two-dimensional cell average in a third-order accurate fashion as well. 
	Naturally we can evaluate the flux integrals using a composite Simpsons method in space and time yielding the following nine-point integral expressions: 
	\begin{align}
		\label{eqn:fluxintegral_stepX}
		G_{i+1/2, j} &\approx \frac{1}{36} (v_{j-1/2}f_{i+1/2,j-1/2}^{*n} + 4v_{j}f_{i+1/2,j}^{*n} + v_{j+1/2}f_{i+1/2,j+1/2}^{*n} \\
		&+ 4v_{j-1/2}f_{i+1/2,j-1/2}^{*n+1/2} + 16v_{j}f_{i+1/2,j}^{*n+1/2} + 4v_{j+1/2}f_{i+1/2,j+1/2}^{*n+1/2} \notag  \\
		&+ v_{j-1/2}f_{i+1/2,j-1/2}^{*n+1} + v_{j}f_{i+1/2,j}^{*n+1} + v_{j+1/2}f_{i+1/2,j+1/2}^{*n+1} ) \notag
	\end{align}
	and 
	\begin{align}
		\label{eqn:fluxintegral_stepV}
		H_{i, j+1/2} &\approx \frac{1}{36} (E_{i-1/2}^* f_{i-1/2,j+1/2}^{*n} + 4 E_{i}^* f_{i,j+1/2}^{*n} + E_{i+1/2}^* f_{i+1/2,j+1/2}^{*n} \\
		&+ 4E_{i-1/2}^{*} f_{i-1/2,j+1/2}^{*n+1/2} + 16 E_{i}^{*} f_{i,j+1/2}^{*n+1/2} + 4 E_{i+1/2}^{*} f_{i+1/2,j+1/2}^{*n+1/2} \notag \\
		&+ E_{i-1/2}^{*} f_{i-1/2,j+1/2}^{*n+1} + 4 E_{i}^{*} f_{i,j+1/2}^{*n+1} + E_{i+1/2}^{*} f_{i+1/2,j+1/2}^{*n+1} \notag
		) \; .
	\end{align} Note that, due to the splitting, the time evolution of the distribution function $f$ only contains the $x$ or $v$ directional advection in \eqref{eqn:fluxintegral_stepX} and \eqref{eqn:fluxintegral_stepV}, respectively. We therefore use the asterisk-superscript in these formulas. In the splitting, we use in each step as $f^{*n}$ the solution from the previous time step or the solution after other split steps that are performed before the current split step. Similarly, $f^{*n+1}$ refers to the solution after this split step which corresponds to $f^{*n}$ in the following split step or to $f^{n+1}$ for the last split-step.
	The central point values $f^{*n}_{i,j+1/2}$ and $f^{*n}_{i+1/2,j}$ in the above integral expression again need to be determined beforehand by evaluating the cell reconstruction at the central position
	  \begin{align}
	  	f^{*n}_{i,j+1/2} = \frac{1}{4} (6\Hat{f}^{*n}_{i,j+1/2} - f^{*n}_{i-1/2,j+1/2} - f^{*n}_{i+1/2,j+1/2}) \; , \\
	  	f^{*n}_{i+1/2,j} = \frac{1}{4} (6\Hat{f}^{*n}_{i+1/2,j} - f^{*n}_{i+1/2,j-1/2} - f^{*n}_{i+1/2,j+1/2}) \; .
	  \end{align}
  	When carrying out only the one-dimensional directional update in the $v$-direction we can recognize that the electric field effectively does not change during that update. 
  	This is due to the fact that the plasma charge density is in fact an integral over velocity space of the distribution function and hence only a function depending on $x$. This means that for our directional splitting approach we can generally expect 
  	a constant electric field during the update step in velocity direction,  simplifying the computation of the horizontal numerical flux (\ref{eqn:fluxintegral_stepV}) by only requiring the solution of the Poisson equation with the current time-step information of $t=t^n$. \\
  	In Alg.~\ref{alg:ThirdOrderFluxIntegral}, find the definition of the operators $L_X$ and $L_V$ for the third-order method. Note that the two-dimensional integrals in \eqref{eqn:fluxintegral_stepX} and \eqref{eqn:fluxintegral_stepV} can be efficiently computed based on sum-factorization to avoid an exponential scaling in the dimension hence the electric field is denoted by $E^*$ at each time level in \eqref{eqn:fluxintegral_stepV}. Whether it equals $E^n$, $E^{n+1}$ or some intermediate value depends on the position of the corresponding split step. 
  	
	\begin{figure}[ht]
		\centering
		\includegraphics[width = 0.6\linewidth]{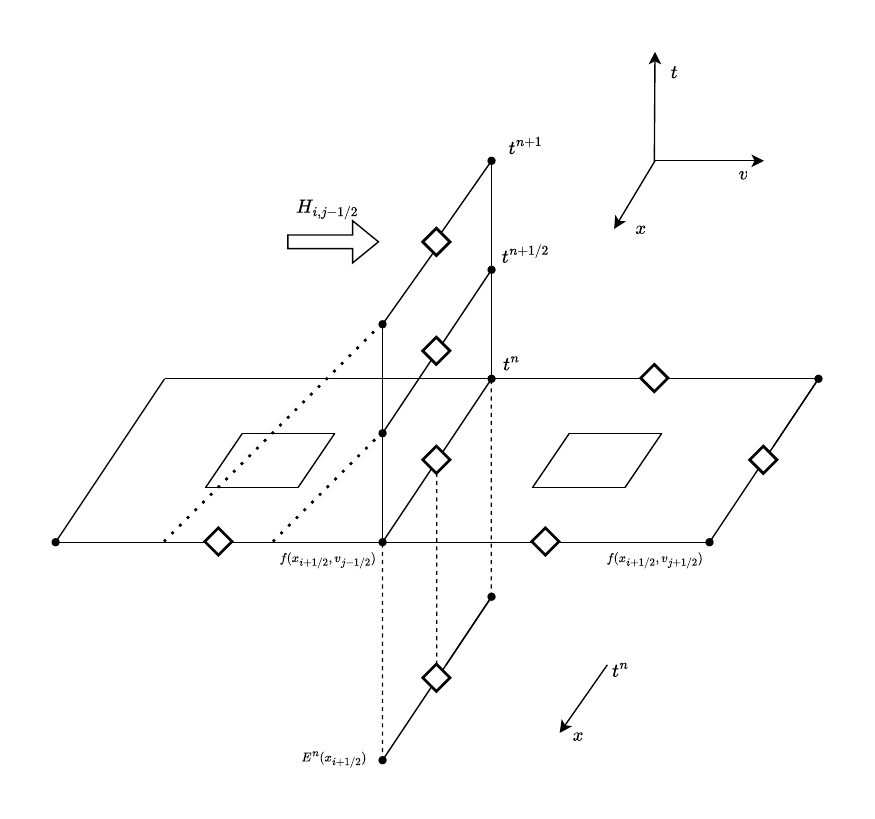}
		\caption{Schematic update procedure for the two-dimensional cell-average. The values located on the $x$-edge at $v_{j-1/2}$ are updated by one-dimensional tracing along the $v$-direction to the time levels $t^{n+1/2}$ and $t^{n+1}$. The interface flux $H_{i, j-1/2}$ for that same edge is obtained by integration over the $x$-$t$-plane using the information of the electric field at time $t^n$.}
		\label{fig:2Dfluxintegral}
	\end{figure}

	\begin{minipage}{0.88\linewidth}
		\newpage
		\begin{algorithm}[H]
			\caption{Third-Order Flux Integral Advection Operators}
			\label{alg:ThirdOrderFluxIntegral}
			\DontPrintSemicolon
			\SetKwFunction{OStepX}{$L_X$}
			\SetKwFunction{OStepV}{$L_V$}
			\SetKwProg{Fn}{Operator}{:}{}
			\Fn{\OStepX{$f$,$\Delta t$}}{
				propagate $f$ over $\Delta t$ / solve: $\partial_t f+ v\partial_x f = 0$\; 
				\begin{itemize}
					\item Apply the update formulae (\ref{eqn:UpdateFormula_1D_Interface}) and (\ref{eqn:UpdateFormula_1D_Average}) to update one-dimensional slices of the grid on the horizontal edges to the time stages $t=t^{n+1/2}$ and $t^{n+1}$ with the nodal values $f_{i+1/2,j+1/2}$ acting as interfaces and the $x$-averages $\hat{f}_{i,j+1/2}$ acting as cell-averages during the update procedure
					\item Apply the update formula (\ref{eqn:UpdateFormula_1D_Interface}) to perform the tracing update of the vertical line-averages $\hat{f}_{i+1/2,j}$ to the time stages $t=t^{n+1/2}$ and $t^{n+1}$
					\item Compute the nine-point flux integral $G_{i+1/2,j}$ (\ref{eqn:fluxintegral_stepX})
					\item Evolve the two-dimensional cell-average according to the directional Finite Volume update
					\begin{align}
						\Bar{f}_{i,j}^{n+1} = \Bar{f}_{i,j}^{n} - \frac{\Delta t}{\Delta x} (G_{i+1/2,j} - G_{i-1/2,j})
					\end{align}    
				\end{itemize}
				calculate density: $\rho(x,t^{n+1}) = \int  f(x,v,t^{n+1}) dv$ \;
				solve Poisson equation: $\partial_x^2 \phi(x,t) = -\rho(x,v,t^{n+1})/\epsilon_0 \; , \;\; E(x,t) = -\partial_x \phi$ \; 
				\begin{itemize}
					\item Apply the scheme from Sec.~\ref{sec:poisson} to compute the Electric field from the distribution function on the inhomogeneous grid. 
					\item The electric field is consequently in the format: point values on the interfaces $E_{i+1/2}$ and $x$-averages $\hat{E}_i$ at the center locations.
				\end{itemize}
				return f, $E$\;
			}
			\Fn{\OStepV{$f$, $E$, $\Delta t$}}{
				propagate $f$ over $\Delta t$ / solve: $\partial_t f + (q_e / m_e) E \partial_v f = 0$\;
				\begin{itemize}
					\item Apply the update formulae (\ref{eqn:UpdateFormula_1D_Interface}) and (\ref{eqn:UpdateFormula_1D_Average}) to update one-dimensional slices of the grid on the vertical edges to the time stages $t=t^{n+1/2}$ and $t^{n+1}$ with the nodal values $f_{i+1/2,j+1/2}$ acting as interfaces and the $v$-averages $\hat{f}_{i+1/2,j}$ acting as cell-averages during the update procedure
					\item Apply the update formula (\ref{eqn:UpdateFormula_1D_Interface}) to perform the tracing update of the vertical line-averages $\hat{f}_{i,j+1/2}$ to the time stages $t=t^{n+1/2}$ and $t^{n+1}$
					\item Compute the nine-point flux integral $H_{i,j+1/2}$ (\ref{eqn:fluxintegral_stepV})
					\item Evolve the two-dimensional cell-average $\Bar{f}_{i,j}$ according to the directional Finite Volume update to the time stage $t=t^{n+1}$
					\begin{align}
						\Bar{f}_{i,j}^{n+1} = \Bar{f}_{i,j}^{n} - \frac{\Delta t}{\Delta v} (H_{i,j+1/2} - H_{i,j-1/2})
					\end{align}     
				\end{itemize}
				return f\;
			}
			\;
		\end{algorithm}
	\end{minipage}
	
	\subsection{Discrepancy Distribution}
	The discrepancy distribution formulation of the Active Flux method yields another scheme that can be straightforwardly applied to higher dimensional cases. In this case the grid is completely homogeneous and made up of point values that are corrected after every directional update such that the correct average is retrieved after integration over the cell domain. 
	Directional updates are carried out as depicted in Fig.~\ref{fig:2nd_order_naive} but now on the homogeneous grid from Fig.~\ref{fig:AF_arrangement}b that consists of corrected point values.  
	Cell averages can be retrieved after every time-step by two-dimensional Simpson integration (\ref{eqn:SimpsonAverage2D}) over the cell domain. \\  
	The entire update procedure to evolve the distribution function in time again consists of one-dimensional directional operations that can be implemented as simple slice-wise operations and hence expanded to higher dimensions easily. Compared to the second-order approach the  computation of electromagnetic fields is simplified because of the homogeneous structure of the grid; since there are only pointlike values in the grid the consistent computation of the moments of the distribution such as the density can be performed by simple integration over all available pointlike DOF along the velocity direction. The resulting charge density is consequently qualitatively made up from pointlike quantities which enables for solving the Poisson equation by nodal Poisson solvers without additional reconstruction operations. 
	Find a algorithmic description of the update procedure in Alg.~\ref{alg:DiscrepancyDistribution}. \\
	\begin{minipage}{0.85\linewidth}
	\begin{algorithm}[H]
		\caption{Discrepancy Distribution Advection Operators}
		\label{alg:DiscrepancyDistribution}
		\DontPrintSemicolon
		\SetKwFunction{OStepX}{$L_X$}
		\SetKwFunction{OStepV}{$L_V$}
		\SetKwProg{Fn}{Operator}{:}{}
		\Fn{\OStepX{$f$,$\Delta t$}}{
			propagate $f$ over $\Delta t$ / solve: $\partial_t f+ v\partial_x f = 0$\; 
			\begin{itemize}
				\item Apply the formulae (\ref{eqn:DiscrepancyFormulaInterfacePosVel})-(\ref{eqn:DiscrepancyFormulaCenterNegVel}) to update horizontal slices of the homogeneous grid  
			\end{itemize}
			calculate density: $\rho(x,t^{n+1}) = \int  f(x,v,t^{n+1}) dv$ \;
			solve Poisson equation: $\partial_x^2 \phi(x,t) = -\rho(x,v,t^{n+1})/\epsilon_0 \; , \;\; E(x,t) = -\partial_x \phi$ \; 
			\begin{itemize}
				\item Apply the scheme from Sec.~\ref{sec:poisson} to compute the electric field from the distribution function on the homogeneous grid. 
				\item The electric field is consequently in the format: point values on the interfaces $E_{i+1/2}$ and in the center $E_i$
			\end{itemize}
			return f, $E$\;
		}
		\Fn{\OStepV{$f$, $E$, $\Delta t$}}{
			propagate $f$ over $\Delta t$ / solve: $\partial_t f + (q_e / m_e) E \partial_v f = 0$\;
			\begin{itemize}
				\item Apply the formulae (\ref{eqn:DiscrepancyFormulaInterfacePosVel})-(\ref{eqn:DiscrepancyFormulaCenterNegVel}) to update vertical slices of the homogeneous grid
			\end{itemize}
			return f\;
		}
		\;
	\end{algorithm}
	\end{minipage}
	
	\section{Consistently Solving the Poisson Equation on Active Flux Grids} \label{sec:poisson}
	To compute the charge density from the plasma distribution function the moment integral in the velocity direction needs to be evaluated. The distribution function on the split step Active Flux grids is inhomogenously made up from point-values, line- and cell-averages, but we require the charge density as point-values-only to afterward consistently apply nodal Poisson solvers such as sufficient order Finite-Difference or Fourier spectral methods for periodic domains. \\ 
	To overcome this inconvenience the moment integral is calculated in the following way:
	\begin{enumerate}
		\item Along the $v$-direction integrate over the vertical line averages that are in fact point values along $x$ to obtain the charge density at the interface positions $x_{i+1/2}$;  
		\begin{align}
			\rho_{i+1/2}^n = \Delta v\sum_{j=1}^{N_v} \Hat{f}_{i+1/2, j}^n \; . 
		\end{align}
		\item Along the $v$-direction integrate over the cell averages $\Bar{f}_{ij}$ to obtain the $x$-averaged charge density
		\begin{align}
			\Hat{\rho}_i^n = \Delta v\sum_{j=1}^{N_v} \Bar{f}_{i, j}^n \; . 
		\end{align}
		Here $\Hat{\rho}_i^n$ denotes the one-dimensional line averaged charge density of the ith cell and $\Bar{f}_{i, j}^n$ . 
		\item Use the one-dimensional reconstruction polynomial (\ref{eqn:Reconstruction_1D}) to obtain $\rho_i^n$ as point values at the cell center by evaluating at $\xi = 0.5$;  
		\begin{align}
			\rho_i^n = \frac{1}{4} (6 \Hat{\rho}_i^n - \rho_{i-1/2}^n - \rho_{i+1/2}^n) \; . 
		\end{align}
	\end{enumerate}
	As a result the charge density $\rho_i$ is of point values only and can be used as a source term for a nodal Poisson solver of choice. 
	The resulting discretized electric field $E_i$ consists also solely of point values and accordingly needs to be averaged to yield the correct shape of point values at the interfaces and line averages along the $x$-direction at the center positions to be applicable in the $v$-step update on the inhomogeneous Active Flux grid 
	\begin{align}
		\Hat{E}_i^n = \frac{1}{6} (E_{i-1/2}^n + 4E_i^n + E_{i+1/2}^n) \; . 
	\end{align}
	For the discrepancy distribution the grid is homogeneously made up of point values; hence we can compute the integrals in $v$-space by applying Simpsons rule over each cell. In this regard the discrepancy distribution method simplifies the computation while still being consistent. Also, this approach is straightforward expandable to higher-dimensional scenarios.
	
	\section{Numerical Experiments}\label{sec:applications}
	In the following section we will test the performance of the previously presented schemes in the context of numerical experiments. For the Poisson equation, we use a Fourier spectral solver in all our experiments. As benchmark problems the following systems from electrostatic plasma physics will be considered in 1D1V phase space: 
	\begin{itemize}
		\item \textbf{Landau Damping:}
		We consider a two-dimensional phase-space domain with $-2\pi \leq x \leq 2\pi$ and $-5 \leq v \leq 5$ where the spatial, velocity and time dimension is normalized with respect to the Debye length, the electron thermal velocity and the electron plasma frequency, respectively. 
		Periodic boundary conditions are imposed on the physical and the velocity space. 
		The initial condition for the electron distribution function is given by a Maxwellian velocity distribution with a small scale perturbation in physical space, 
		\begin{align}
			f_e(x, v, t=0) = \frac{1}{\sqrt{2 \pi}} \exp{\left(-\frac{v^2}{2} \right)}  (1 + A \cos(kx)) .
			\label{eq:LD_initial_condition}
		\end{align}
		We consider two choices of the parameters. For the \textbf{Weak Landau damping} test case, the amplitude and the wavenumber of the perturbation are $A = 10^{-3}$ and $k = 0.5$. The CFL number is for all simulations set to $\nu = 1/\pi \approx 0.318 < 0.5$ to also satisfy the stricter time-step restriction of the discrepancy distribution method. \\
		For this physical system analytical theory predicts a characteristic exponential damping $\propto \exp(- \gamma t)$ and oscillation of the electric fields energy with a decrement of $\gamma \approx 0.153$. \\
		By increasing the amplitude of the initial perturbation for the weak Landau damping problem (\ref{eq:LD_initial_condition}) the system develops a non-linear behavior in the sense that there is no constant damping rate for the electric energy. 
		This problem will be denoted in the following as \textbf{Strong Landau Damping}.
		For this problem again the initial condition (\ref{eq:LD_initial_condition}) is considered over the same periodic computational domain but with $A = 0.5$. Other parameters such as the time-step are chosen the same as for the previous weak Landau damping test case. 
		\item \textbf{Two Stream Instability:} 
		Furthermore the unstable system of two counterstreaming beams in a one-dimensional plasma with a small initial perturbation is considered. The electric field energy is expected to increase exponentially over time. The initial electron distribution function is set up to be 
		\begin{align}
			f_e(x, v, t=0) = \frac{1}{2\sqrt{2 \pi}} \left[ \exp{\left( \frac{-(v-v_0)^2}{2} \right)} + \exp{ \left( \frac{-(v+v_0)^2}{2} \right)} \right] (1 + A\cos(kx))
		\end{align}
		over the periodic phase space volume $-5\pi \leq x \leq 5 \pi$ and $-10 \leq v \leq 10$. The beam velocities are set to $v_0 = 3$ and the perturbation has a wavenumber of $k=0.2$ and an amplitude of $A=10^{-3}$.
	\end{itemize}
	We furthermore test against a split-step Vlasov--Poisson solver based on the \emph{Positive and Flux Conservative (PFC)} method \cite{filbet-sonnendruecker-bertrand:2001} as a state-of-the-art numerical scheme. The PFC method is essentially a third-order accurate semi-Lagrangian method that preserves positivity of the distribution function. The computational grid for PFC consists solely of cell averages that are evolved in time; no additional cell-DOF are introduced which reduces the effective numerical cost as well as the memory footprint in comparison to the Active Flux methods with the same number of cells. Note that $N_x$, $N_v$ refer to the number of \emph{cells} in $x$ and $v$ direction, respectively. Hence, the number of DOF is $2N_x$ and $2N_v$ for the Active Flux schemes while it is only $N_x$ and $N_v$ in PFC.\\
	
	\subsection{Convergence}
	We verify the general order of accuracy by simulating two-dimensional classical benchmark problems from electrostatic plasma physics namely weak Landau damping and the Two Stream Instability. \\
	As an error measurement we compute the relative $L_1$-Norm error at a fixed time $t_{max}$. Since there are no analytical solutions available as references we conduct the error computation in the following way:
	Firstly, we perform a high resolution simulation $N_{x,\text{ref}} =  N_{v,\text{ref}} = 512 = 2^9$ of the system with a fixed CFL number up to a given time $t_{\max}$ and take this to be our numerical reference solution $f_{\text{ref}}$. To effectively compute an error between the reference and a lower-resolution simulation result we scale the reference solution down to the desired lower resolution $N_x, N_v$ of a run by recursively computing means of the cell averages. The scaled-down version of the reference solution will be in the following denoted by $\Tilde{f}_{i,j}(t)$. E.g. scaling the reference solution to a grid with a resolution of $N_{x,\text{ref}}/2, N_{v,\text{ref}}/2$ can be executed by computing the mean of four neighboring cell averages of the reference solution as 
	\begin{align}
		\Tilde{\Bar{f}}_{i,j} = \frac{1}{4}(\Bar{f}_{ref_{i,j}} + \Bar{f}_{ref_{i+1,j}} + \Bar{f}_{ref_{i,j+1}} + \Bar{f}_{ref_{i+1,j+1}}) \; . 
	\end{align}
	This procedure can be recursively repeated to achieve down-scaling to any given grid resolution $N_x, N_v = 2^k, 2^k$ for $k < 9$. 
	We define the relative $L_1$-Norm error $\epsilon_{VP}$ as 
	\begin{align}
		\epsilon_{VP}(t) = \frac{\sum_{i=1}^{N_x} \sum_{j=1}^{N_v}|\Bar{f}_{i,j}(t) - \Tilde{\Bar{f}}_{i, j}(t)|}{\sum_{i=1}^{N_x} \sum_{j=1}^{N_v} | \Tilde{\Bar{f}}_{i, j}(t)|} \; . 
	\end{align}
	where $\Bar{f}_{i,j}$ are the cell averages of the numerical solution with a given resolution $N_x, N_v$.  
	
\begin{figure}[H]
		\centerline{\includegraphics[scale=0.45]{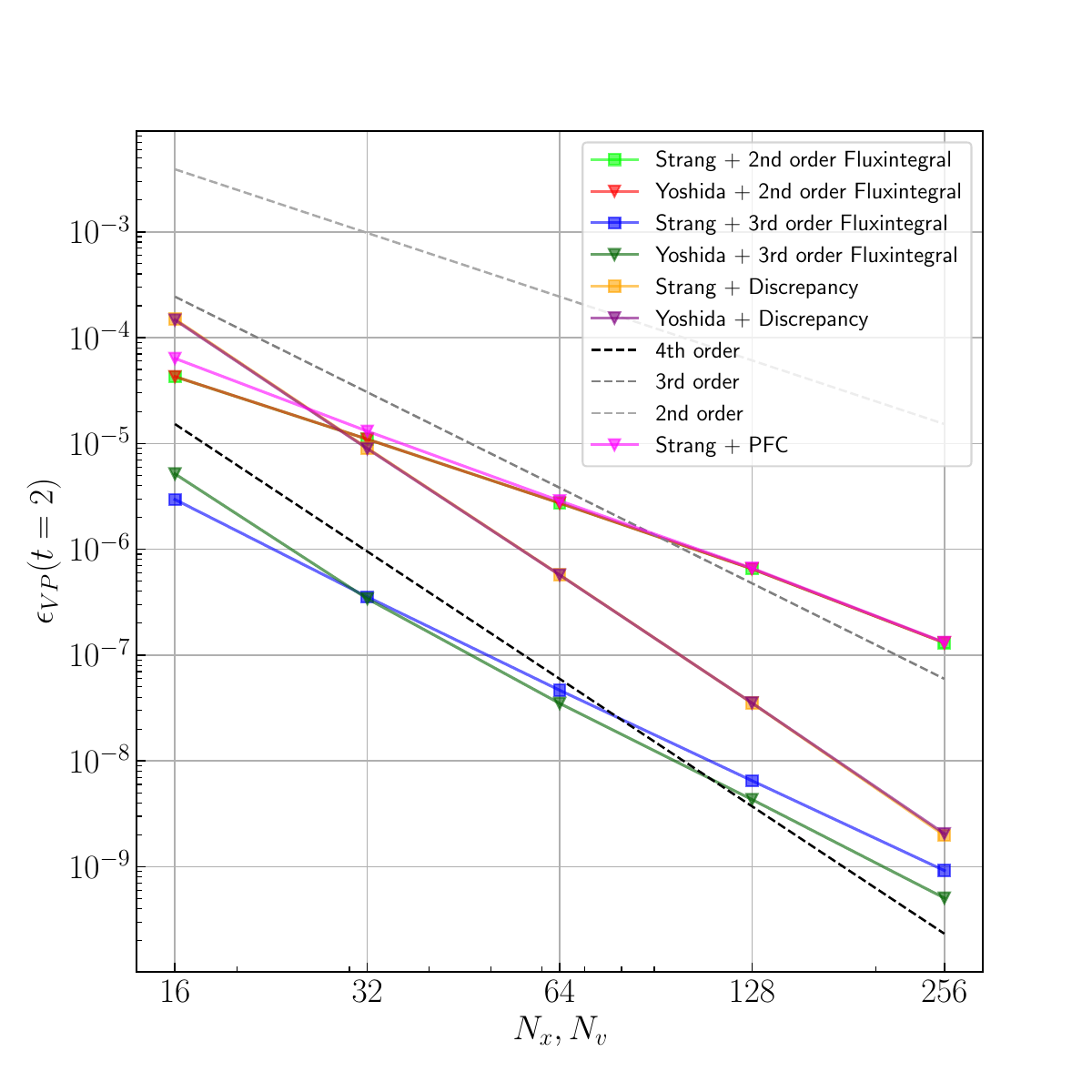}
			\includegraphics[scale=0.45]{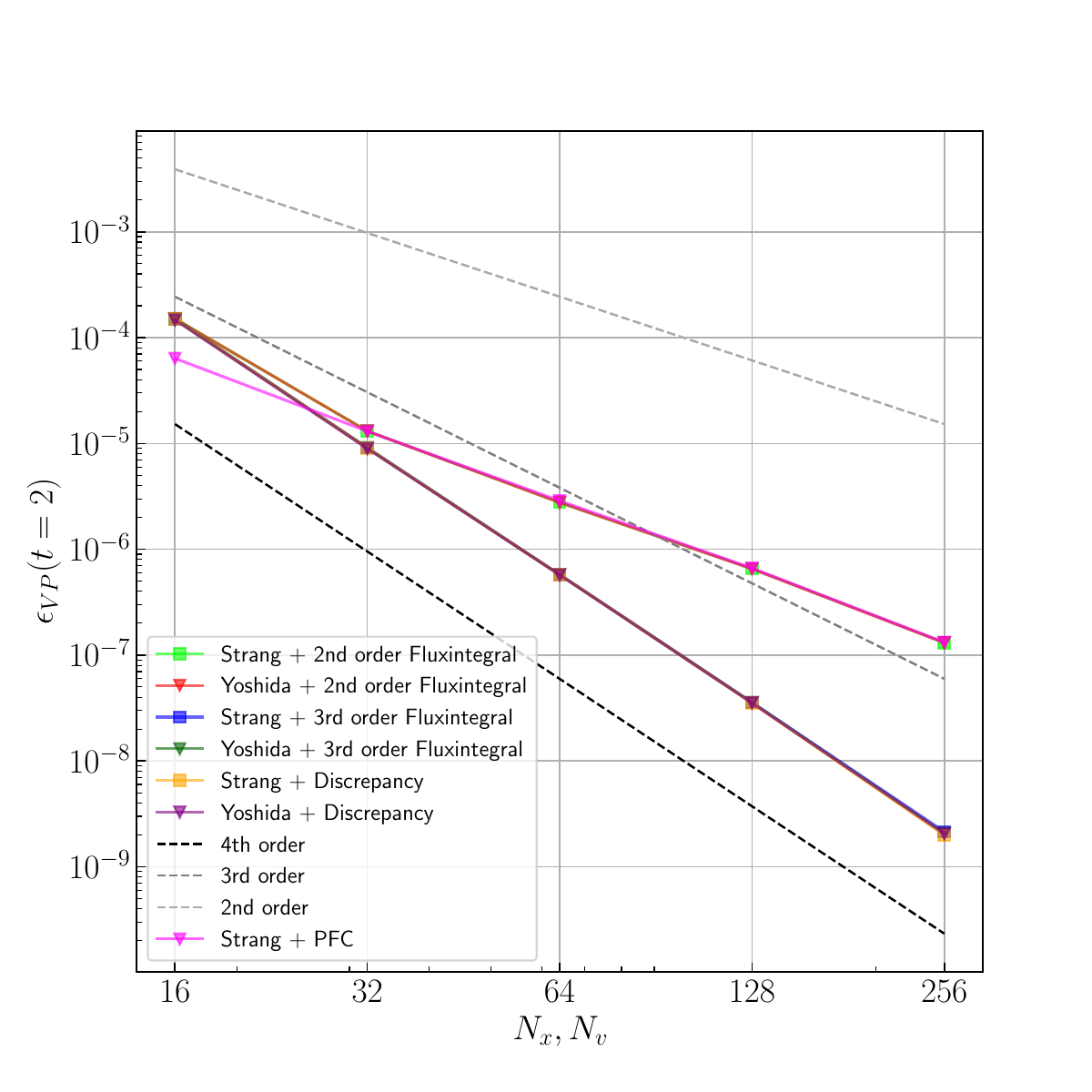} 
		}
		\caption{Convergence of error for the initial phase of the weak Landau damping problem. Integrals were carried out analytically (left) or numerically by Simpsons rule (right) for schemes that require initial averaging of the distribution function.}
		\label{fig:convergence_LD}
	\end{figure}

\begin{figure}[H]
		\centerline{\includegraphics[scale=0.45]{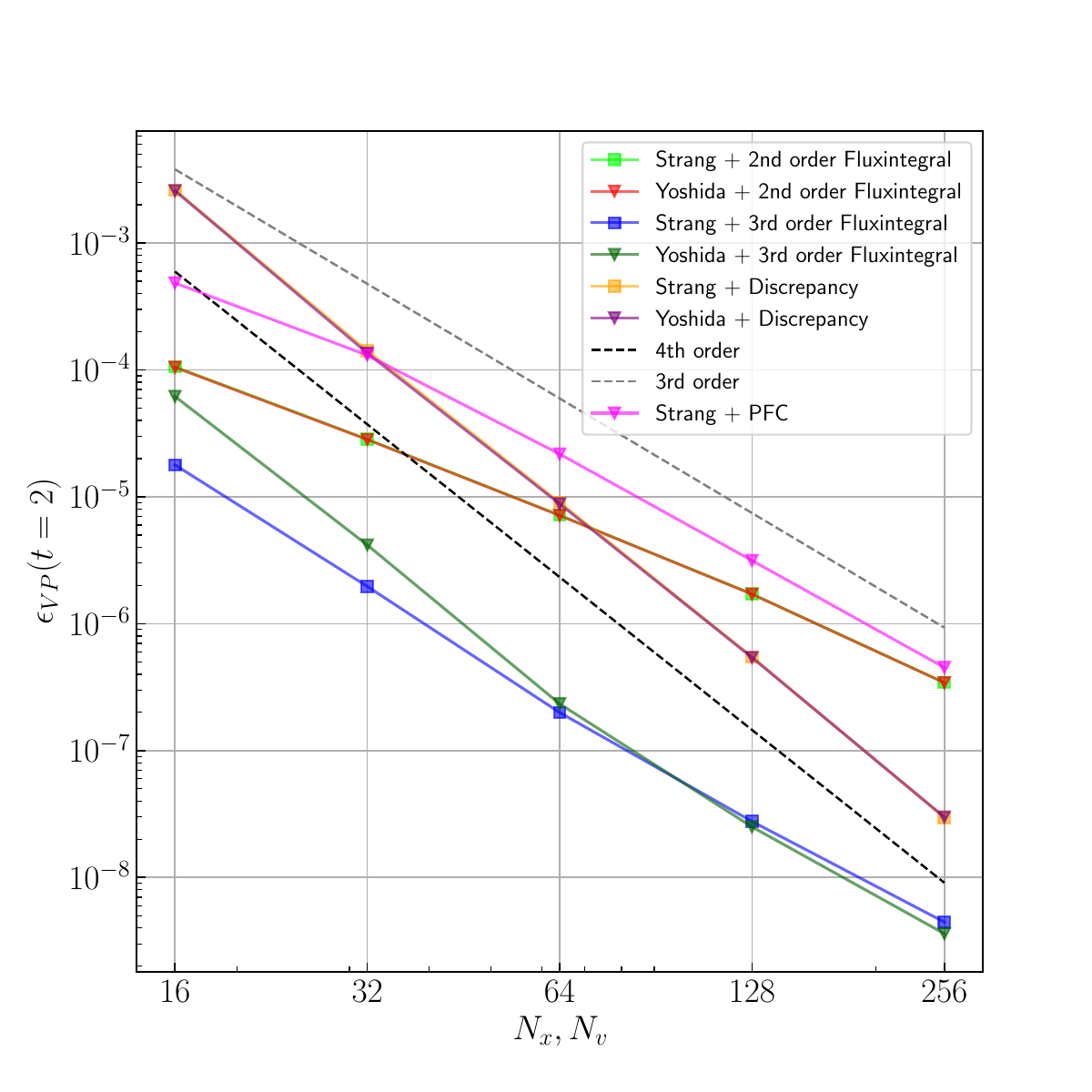}
			\includegraphics[scale=0.45]{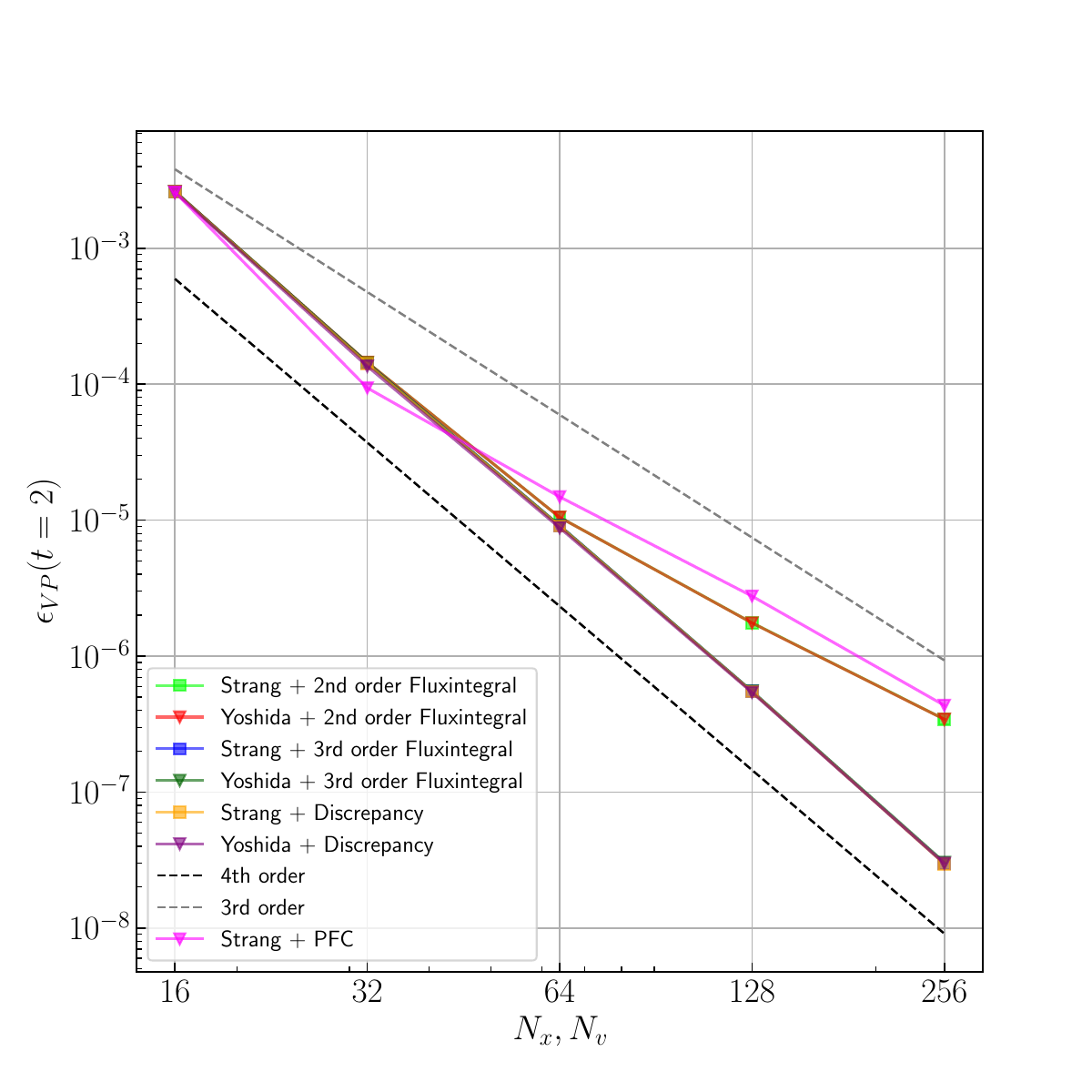} 
		}
		\caption{Convergence of error for the initial phase of the Two-Stream instability. Integrals were carried out analytically (left) or numerically by Simpsons rule (right) for schemes that require initial averaging of the distribution function.}
		\label{fig:convergence_TS}
	\end{figure}

	We compare the convergence of the error at time $t_{\max}= 2$ for second (Alg.~\ref{alg:Strang}) and fourth order (Alg.~\ref{alg:Yoshida}) accurate operator-splitting in time as well as for the various proposed approximations to the numerical flux integrals in space. Results of weak Landau damping are shown in Fig.~\ref{fig:convergence_LD} and for the Two Stream Instability in Fig.~\ref{fig:convergence_TS}. 
	As a reference we additionally plot a curve produced by the PFC method~\cite{filbet-sonnendruecker-bertrand:2001} in combination with Strang operator splitting with the same number of cells for both problems. All integrals that are required to set up the grids for the different methods were computed analytically in this case. We find the expected convergence behavior for the proposed second- and third-order flux integrals with even better than third-order convergence for the discrepancy distribution method in the weak Landau damping test case.  
	The influence of the time splitting scheme can be generally found minor as it only shows a difference for low spatial errors. It can be concluded that the spatial error emerging from the employed advection scheme is the dominant factor in this kind of convergence studies. \\
	It was furthermore observed that the method for computing the integrals in the initial setup of the grid has an influence on the error constant as we are just simulating the initial phase of the problem: In comparison the alternative convergence curves where the initialization of the integrals in the grid has been carried out numerically using Simpsons rule over the one-dimensional edges or two-dimensional cells respectively are shown in Fig.~\ref{fig:convergence_LD}. Due to the influence of the grid initialization the curves for the different advection methods are more similar while laying for most part on top of each other. There is only a slight noticeable branching between the proposed third- and second-order accurate methods with the third-order flux integral approximation showing throughout even better than third-order accuracy in our test. \\
	It must be noted that the discrepancy distribution method does not require the initial computation of integrals for the grid setup as the grid is initialized with point values only. The computation of the discrepancy during runtime is executed using Simpsons rule and hence cannot be executed analytically. For comparability, we display the same discrepancy-distribution curve in both figures for every problem.  \\
	
	\subsection{Long Time Behavior of the Electric Field}
\begin{figure}[ht]
		\centering
		\includegraphics[width = 0.6\textwidth]{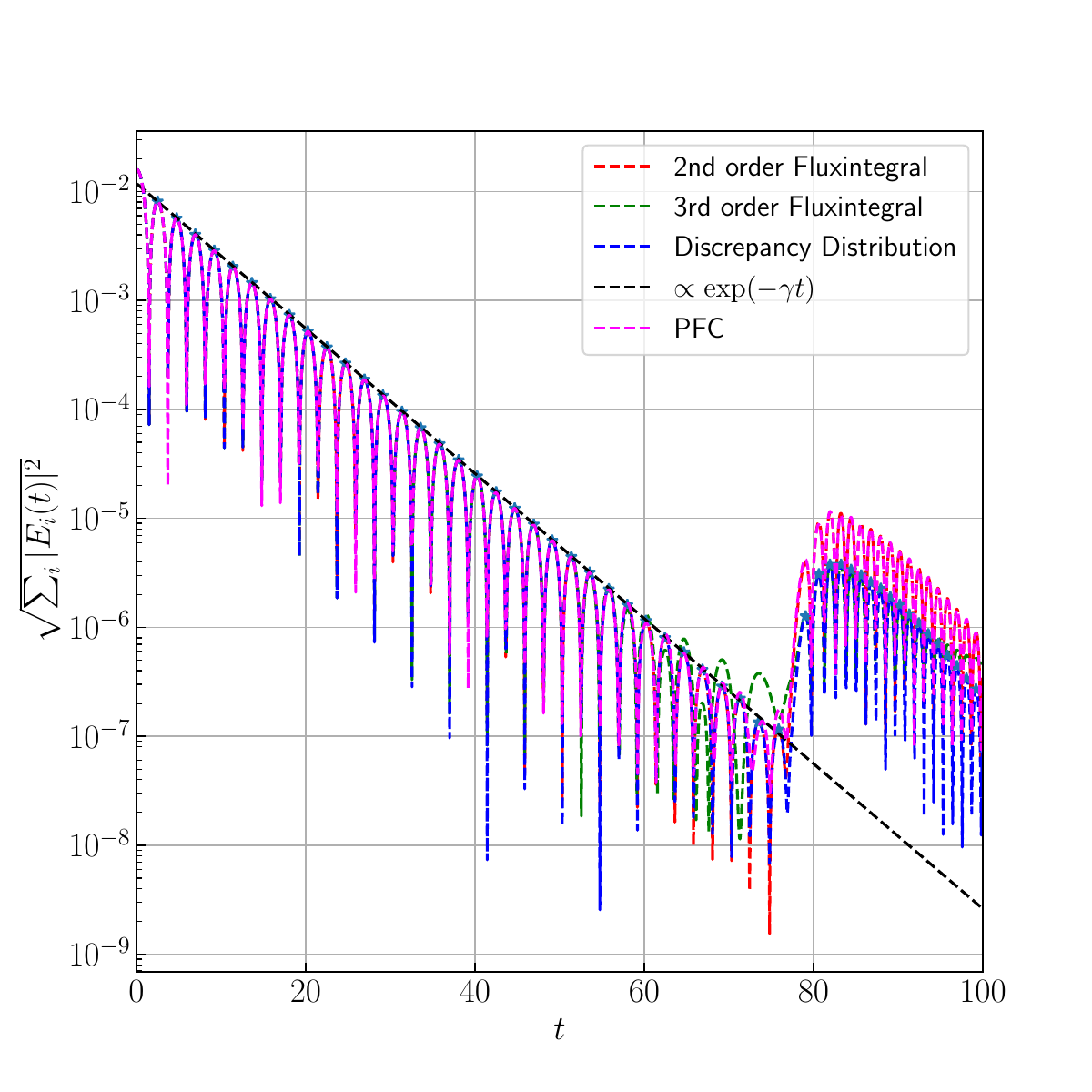}
		\caption{Comparison of weak Landau damping oscillations for different methods, on $N_x, N_v = 128, 128$ grid.  
	}
		\label{fig:LongLD}
	\end{figure}

\begin{figure}
		\centerline{\includegraphics[scale=0.45]{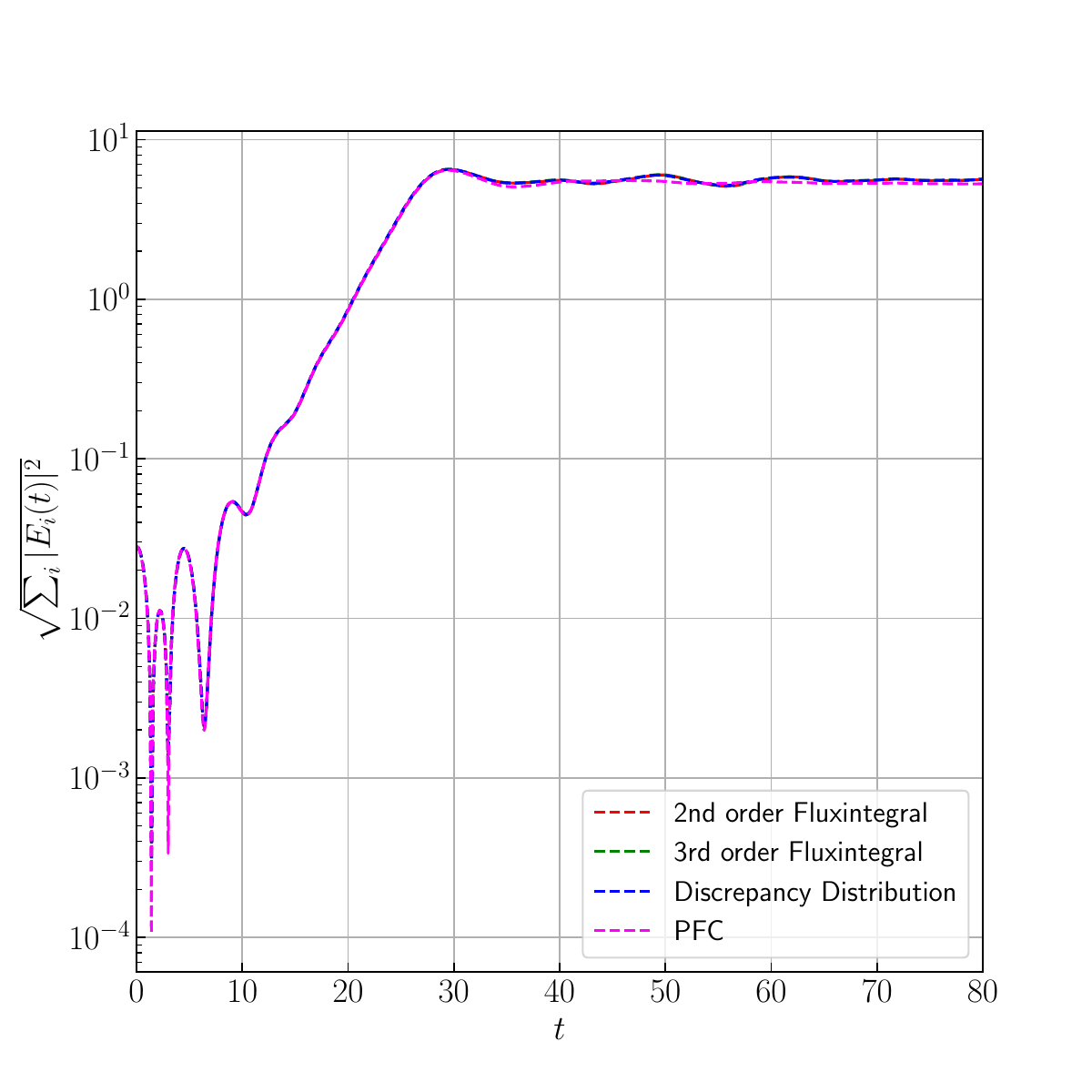}
			\includegraphics[scale=0.45]{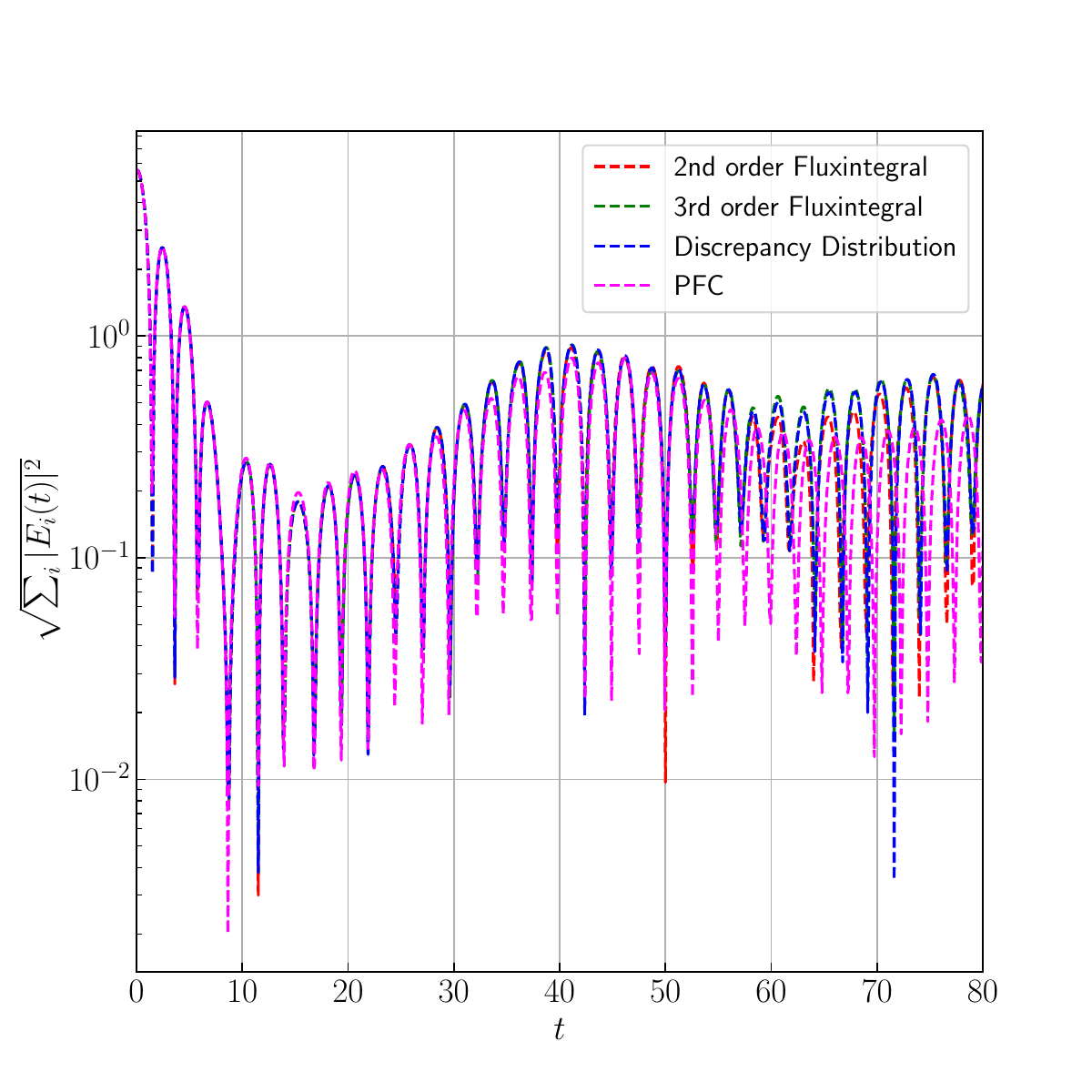}}
		\caption{Evolution of the Electric field energy for the Two Stream Instability (left) and the strong Landau damping (right) problem. Both simulations were carried out with $N_x, N_v = 64, 64$.}
		\label{fig:TS_SLD}
	\end{figure}

	Another aim for a numerical scheme is to achieve long-time stability of the simulation. When performing weak Landau damping simulations the free streaming part of the transport equation drives the distribution function to develop filaments that become increasingly fine over time.  
	Recurrence effects of the electric field start to manifest as a strong exponential increase of the electric field contrary to the analytical prediction once the highly filamented distribution function is no longer sufficiently resolved numerically on the grid. \\
	We plot the electric energy long-time behavior for the weak Landau damping problem simulated by Active Flux and PFC methods and Strang splitting in Fig.~\ref{fig:LongLD}. All the Active Flux methods show good agreement with the analytical damping rate (dashed black line) as well as with the reference PFC-solution. The recurrence times differ with the third-order flux integral method showing slight deviations from the characteristic damping behavior first. 
	The second-order flux integral approximation produces very similar results to the reference PFC method. \\
	Additionally, a similar test has been carried out for the Strong Landau problem in Fig.~\ref{fig:TS_SLD}. The characteristic damping and growth behavior of the electric field energy in the beginning phase of strong Landau damping  is well reproduced by all the methods. At later time the stronger numerical dissipation for the PFC result shows in a significant decrease in the energy compared to the Active Flux methods results.  
	Furthermore, in Fig.~\ref{fig:TS_SLD} the development of the electric field energy for the unstable Two Stream system is shown. The exponential increase of the electric energy is reproduced by all the Active Flux methods showing strong resemblance to the reference PFC solution.
	
	\subsection{Snapshots of the distribution function}
\begin{figure}
		\vspace{-2cm}
		\includegraphics[scale = 0.42]{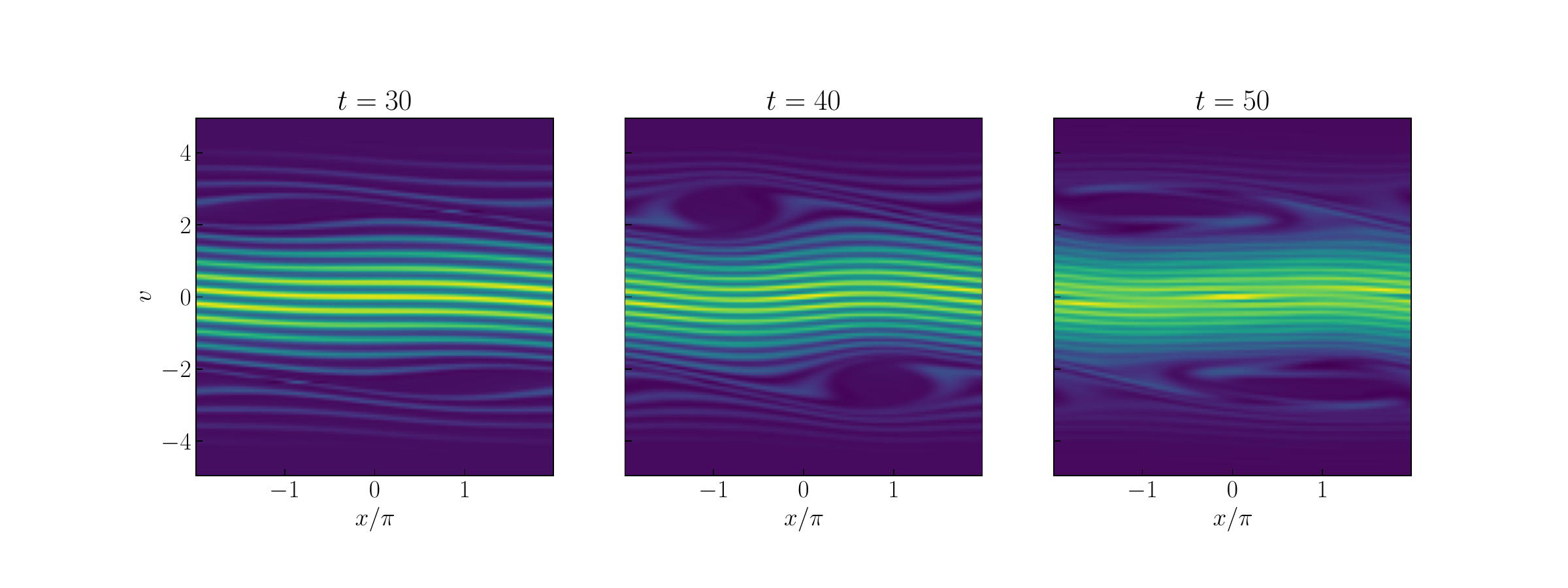}
		
		\vspace{-1cm}
		
		\includegraphics[scale = 0.42]{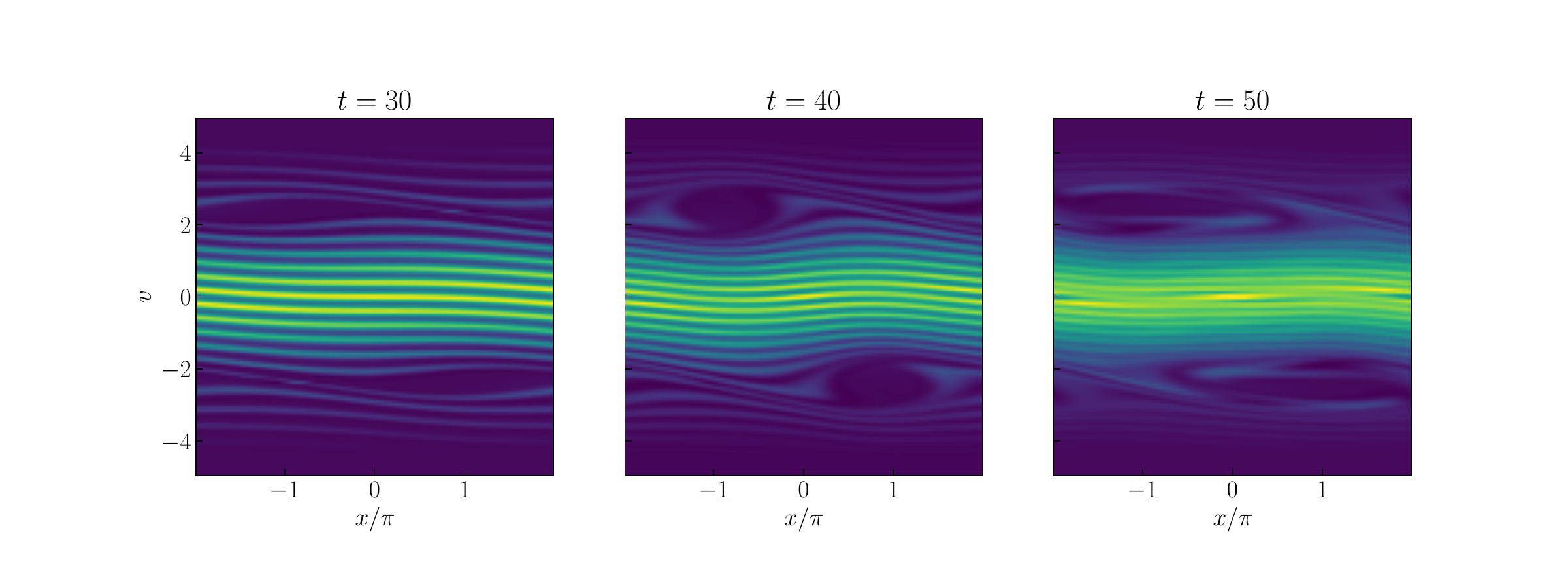}
		
		\vspace{-1cm}
		
		\includegraphics[scale = 0.42]{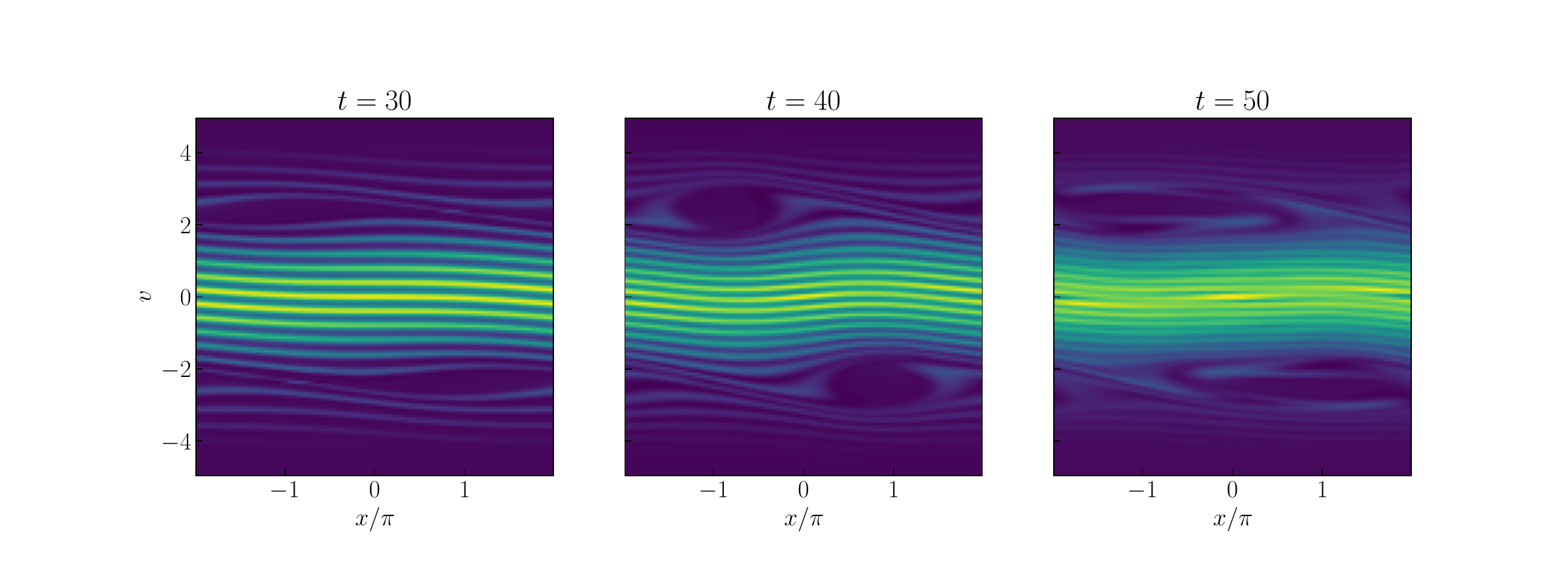}
		
		\vspace{-1cm}
		
		\includegraphics[scale = 0.42]{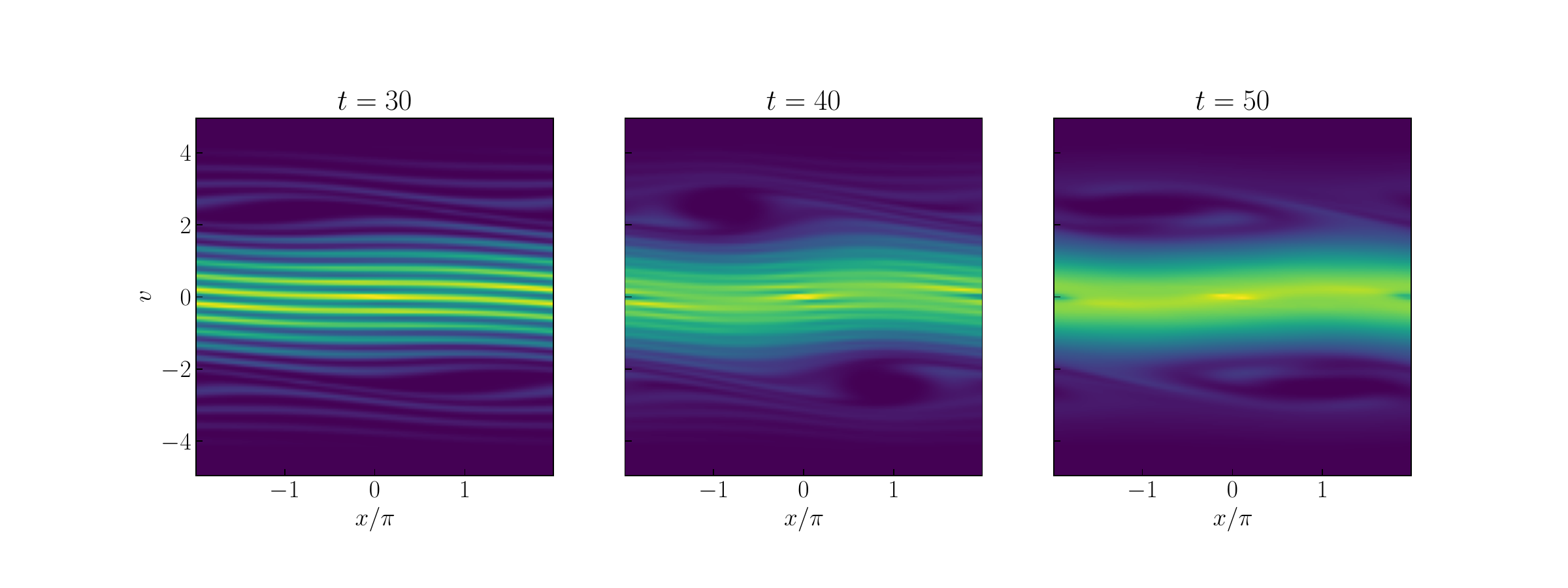}
		
		\caption{Snapshots at different time stages for strong Landau damping. Row-wise from top to bottom: Second-order flux integral AF (first), discrepancy distribution AF (second) and third-order flux integral AF (third) and PFC (fourth) with a resolution of $N_x, N_v = 128, 128$. Only the cell averages are displayed.}
		\label{fig:SLD_snapshots}
	\end{figure}

\begin{figure}[ht]
		\centering
		\includegraphics[scale = 0.42]{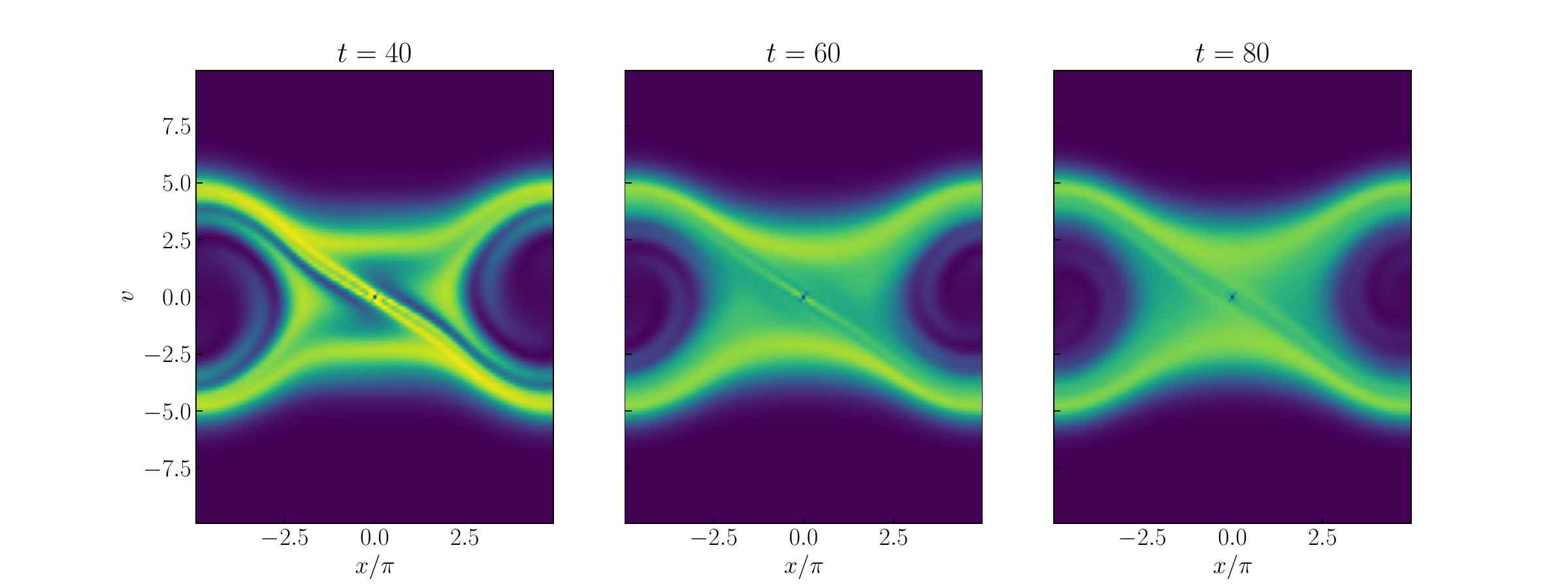}
		
		\vspace{-1cm}
		
		\includegraphics[scale = 0.42]{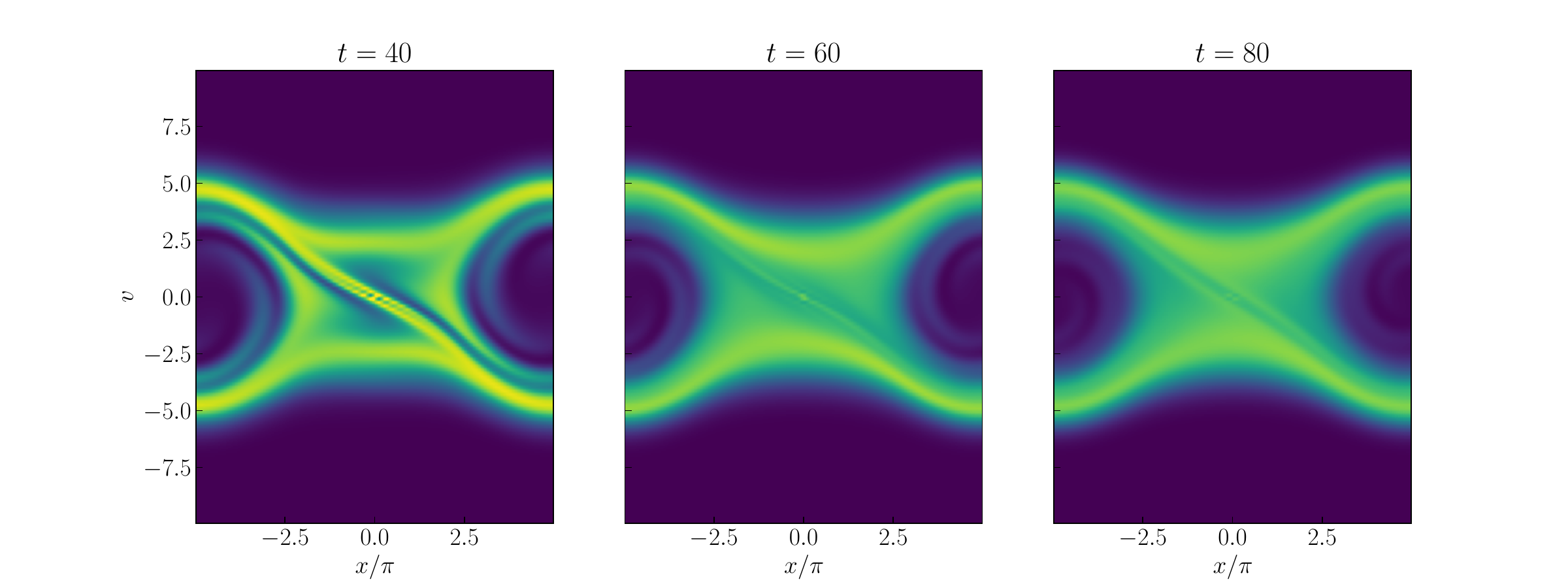}
		\caption{
		Snapshots at different time stages for the Two Stream Instability: Third-order flux integral (top) and PFC method (bottom). Same number of DOF for Active Flux and PFC methods are considered; for AF the original grid resolution was $N_x, N_v = 64, 64$ and $N_x, N_v = 128, 128$ for PFC. The time step for the Active Flux method could be consequentially chosen twice as large as for the PFC method (both with $|\nu| = 1$). 
		PFC cell averages as well as Active Flux line and cell averages are transformed into point values via histopolation.
		}
		\label{fig:TS_snapshots_hist}
	\end{figure}

\begin{figure}[ht]
		\centering
		\includegraphics[scale = 0.35]{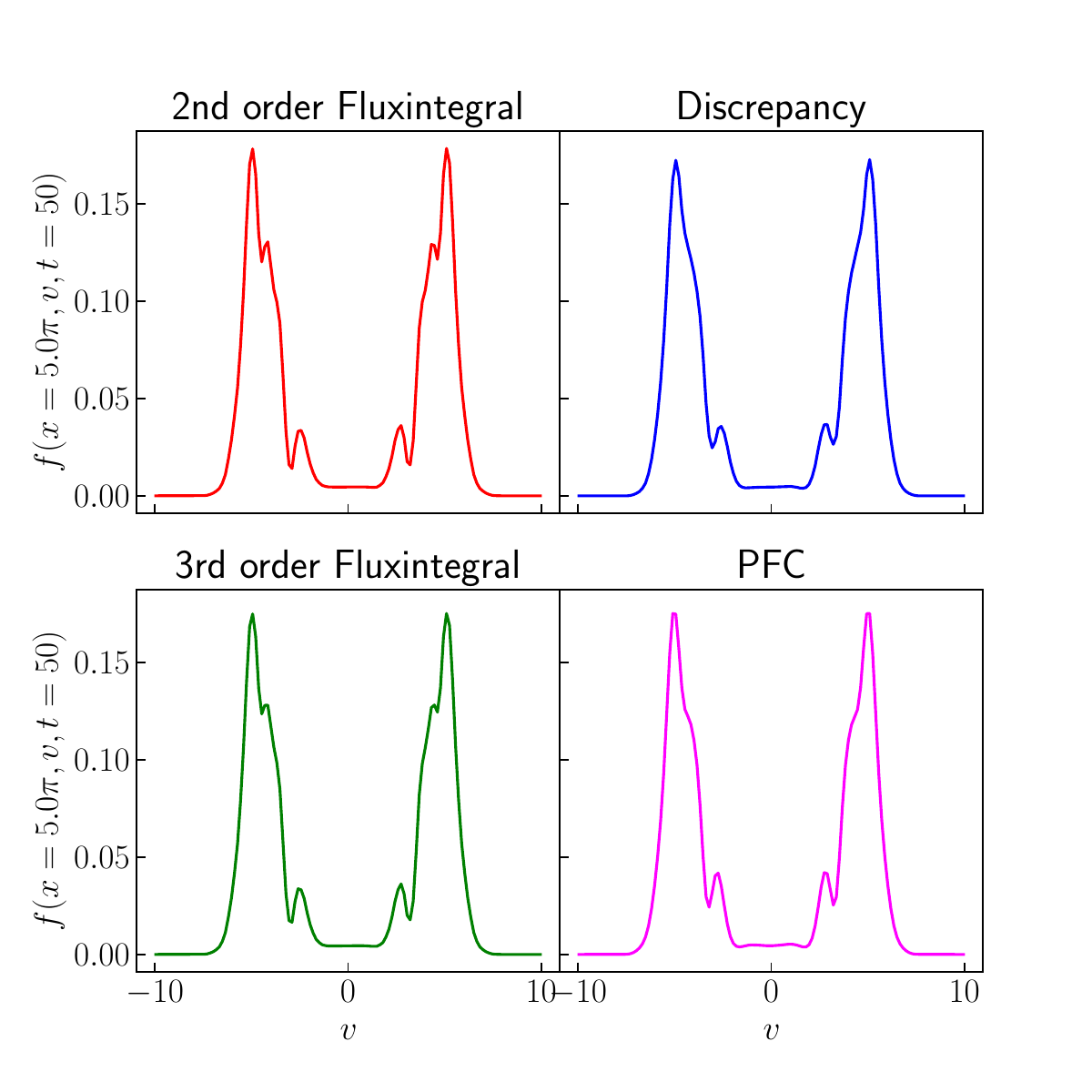}
		\includegraphics[scale = 0.35]{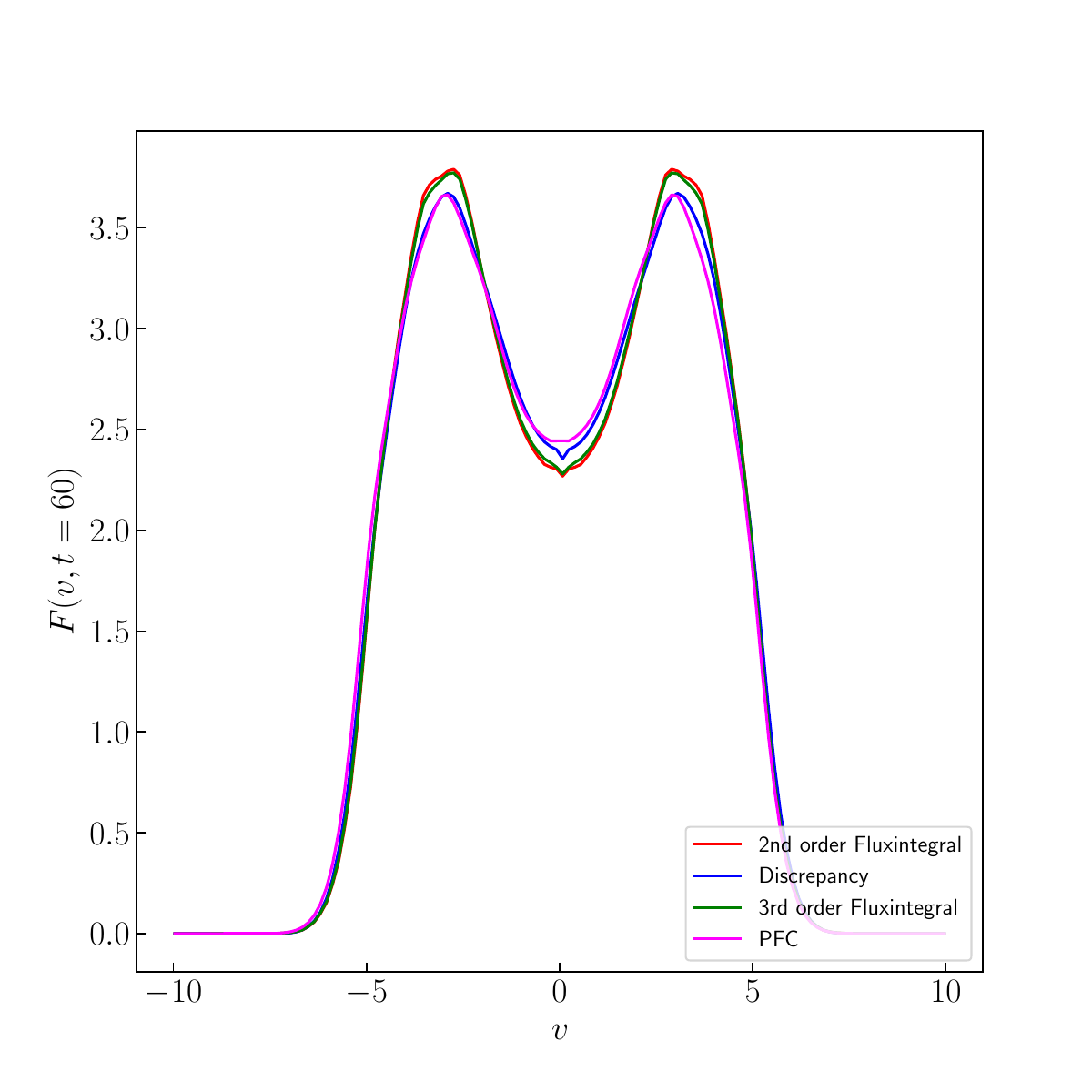}
		\caption{(left) Velocity space cross-section through snapshots at $x=5\pi$ and $t = 50$ for the Two Stream Instability. The grids have the same number of DOF: $N_x, N_v = 64, 64$ for AF and $N_x, N_v = 128, 128$ for PFC. 
		(right) Integrated velocity distribution $F(v,t)$ for the same simulation. }
		\label{fig:cross_sections}
	\end{figure}

	Additionally, we plot snapshots of the distribution function for strong Landau damping at given time stages in Fig.~\ref{fig:SLD_snapshots}. The distribution function develops a fine-scale structure over time that optimally should be reproduced by the numerical method but naturally falls victim to the numerical dissipation to some degree by effectively smearing out the solution over time. In Fig.~\ref{fig:SLD_snapshots} we display the cell averages of the discretized distribution function for every method. One has to note that this comparison is not completely fair, since the Active Flux methods possesses additional DOF per cell which increases the computational-cost in comparison to PFC. \\
	It can be found that all the Active Flux methods produce optically very similar results while being throughout less dissipative at later times than PFC where a lot of the fine scale structure is smeared out due to numerical dissipation. \\
	 Since the AF grids posses additional DOF in the grid that increase memory cost in comparison to the PFC methods grid that only consists of cell averages. Hence, we compare Active Flux simulation runs with $N_x, N_v$ grid cells to a PFC run with $2N_x, 2N_v$ grid cells which results in an equal number of grid DOF. To achieve comparable pictures we use directional fourth order accurate histopolation \cite{robidoux2008polynomial, kormann2024dual} to effectively turn cell averages into corresponding cell centered point values (\ref{eq:2D_histopolation}) for the visualization of PFC results:
	 \begin{align}
	 	f_{i,j} &\approx \frac{1}{576} (676 \Bar{f}_{i,j} - 26 (\Bar{f}_{i+1,j} + \Bar{f}_{i-1,j} + \Bar{f}_{i,j+1} + \Bar{f}_{i,j-1}) + (\Bar{f}_{i+1,j+1} + \Bar{f}_{i+1,j-1} + \Bar{f}_{i-1,j+1} + \Bar{f}_{i-1,j-1}))
	 	\label{eq:2D_histopolation} \\
	 	f_{i,j+1/2} &\approx \frac{1}{24} (26 \hat{f}_{i,j+1/2} - \hat{f}_{i+1,j+1/2} - \hat{f}_{i-1,j+1/2}), \quad f_{i+1/2,j} \approx \frac{1}{24} (26 \hat{f}_{i+1/2,j} - \hat{f}_{i+1/2,j+1} - \hat{f}_{i+1/2,j-1}) \label{eq:1D_histopolation}
	 \end{align}
	 For the Active Flux method we also use histopolation (\ref{eq:2D_histopolation}) for cell averages as well as the one-dimensional histopolation formula (\ref{eq:1D_histopolation}) for line averages. 
	 Furthermore for our comparison we consider the PFC and Active Flux methods at the maximum time step that still guarantees equally expensive periodic boundary exchange of two DOF per update. For the Active Flux method this coincides with the schemes CFL restriction of $|\nu| \leq 1$ or $|\nu| \leq 1/2$ for the discrepancy distribution method. The PFC method is theoretically unrestricted in terms of the time step. However increasing the time step for the semi-Lagrangian PFC method also increases the domain of dependence. In a parallel implementation this results in a larger number of DOF that need to be exchanged on the boundary during every update step. For a six-dimensional simulation the boundary exchange has been found to be the main computational bottleneck (cf.~\cite{schild-raeth-etal:2024}). Consequentially the time step for the third-order PFC method is also restricted by $|\nu| \leq 1$ in this test to achieve equally expensive boundary exchange as for the Active Flux methods. This means that for the considered same-DOF comparison the time step is effectively chosen twice as large for the second- and third-order flux integral as well as equally large for the discrepancy distribution as the one used for PFC.
	The resulting comparative same-DOF snapshots for the Two Stream instability are shown in Fig.~\ref{fig:TS_snapshots_hist}. \\
	We can see that the solution obtained with the Active flux method with third-order flux integral (as well as the other Active Flux solutions, not shown here) show sharper contours in the snapshots. Moreover, as described the second- and third-order flux integral Active Flux methods were obtained for a time step twice as large at equal DOF grids with equally expensive boundary exchange. \\
	We additionally plot velocity cross-sections through the snapshots from simulations with $N_x, N_v = 64, 64$ for Active Flux and $N_x, N_v = 128, 128$ for PFC to compare the level of dissipation in Fig.~\ref{fig:cross_sections}. All methods were again considered at the maximum possible time step.
	It can be observed that the discrepancy distribution method produces a smeared-out solution with less steep gradients compared to the other methods.  
	This is due to the distribution of cell-integrated quantities to point valued DOF to achieve conservation, creating small deviations to the results of the semi-Lagrangian predictor step.
	Furthermore we find less detailed solutions for PFC compared to the second- and third-order flux integral Active Flux methods that operate on inhomogeneous grids.
	In Fig.~\ref{fig:cross_sections} (right) we also find by displaying the integrated velocity distribution $F(v,t) = \int_{\Omega_x} f(x,v,t) dx$ that the discrepancy distribution method possesses similar levels of dissipation as PFC for our test case.      

	\subsection{Conservation}
\begin{figure}[ht]
		\centerline{\includegraphics[scale=0.35]{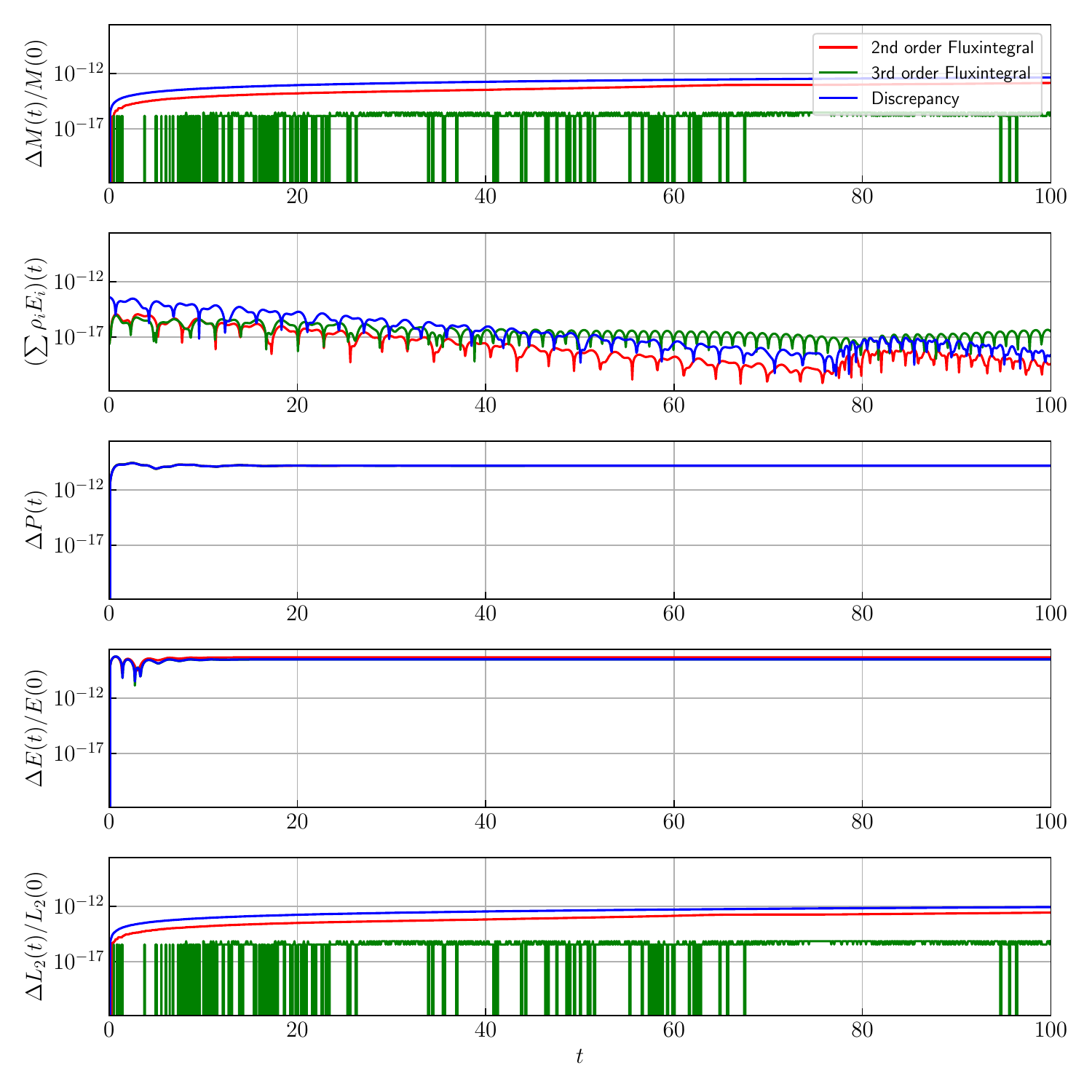}
			\hspace{-1cm}
			\includegraphics[scale=0.35]{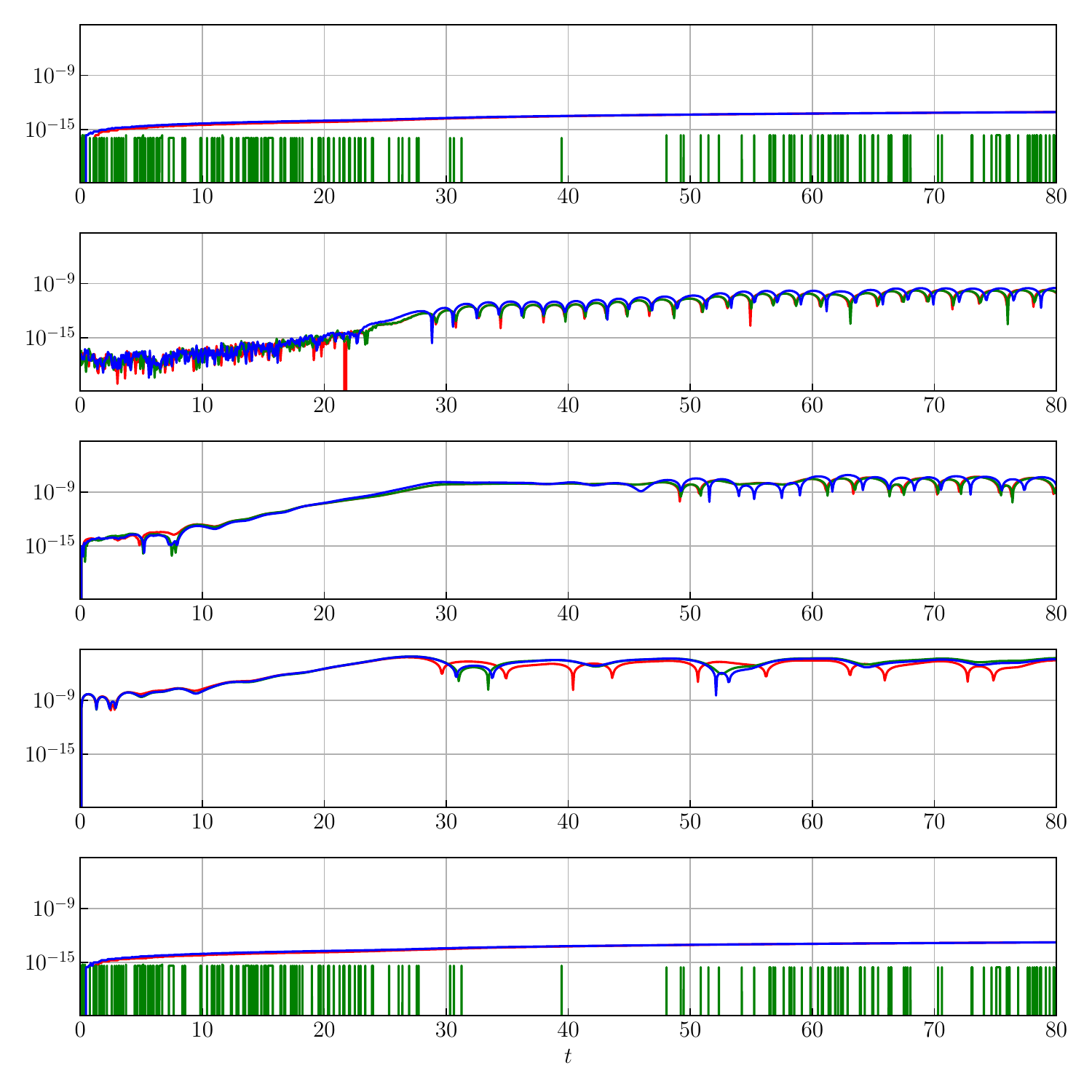}
		}
		\caption{Numerical error of analytically conserved quantities for weak Landau damping (left) and Two-Stream Instability (right) over time. Additionally, the condition $\sum \rho_i E_i = 0$ is a necessary criterion to ensure momentum conservation.}
		\label{fig:Conservation}
	\end{figure}

	One further goal for high-order numerical methods is to replicate the conservative properties of the underlying physical system on a discrete level. For the Vlasov--Poisson system several conserved quantities can be derived from analytical considerations; those are e.g.
	\begin{itemize}
		\item the $L_p$ norms of the distribution function $\int_{\Omega_v} \int_{\Omega_x} |f(x,v,t)|^p dx dv, \quad 1 \leq p < \infty$
		\item the mass $M(t) = \int_{\Omega_v} \int_{\Omega_x} f(x,v,t) dx dv$
		\item the total momentum $P(t) = \int_{\Omega_v} \int_{\Omega_x} v f(x,v,t) dx dv$
		\item  the total energy $W(t) = \int_{\Omega_v} \int_{\Omega_x} v^2 f(x,v,t) dx dv + \int_{\Omega_x} E^2(x, t) dx$. 
	\end{itemize}  
	Furthermore it can be shown that for a split step semi-Lagrangian method the discrete momentum $\Delta x \Delta v\sum_{i,j}^{N_x,N_v} v_j f_{i,j}$ is conserved provided the Poisson solver satisfies $\sum_{i}^{N_x} \rho_i E_i = 0$ with the discrete density $\rho_i = 1 -\Delta v \sum_j \Bar{f}_{i,j}$ (see \cite{sonnendrucker2013numerical}).\\
	In Fig.~\ref{fig:Conservation} the errors over time of some analytically conserved quantities for the weak Landau damping and Two Stream Instability problems are plotted. Exemplary the error of the mass is thereby computed as the absolute difference of the numerical mass and the mass of the initial condition $\Delta M(t) = |M(t) - M(0)|$. For the relative errors we furthermore divide by the initial mass $M(0)$. \\

	In theory, all Active Flux methods conserve both mass and $L_2$ norm up to machine precision.
We note, however, that the use of the update formulas for the update of the cell average leads to a numerical breach of the conservation property because the operations involved in the computation of the flux at a given interface are performed in different order for each adjacent cell due to the factorization of the formulas. This is the case for the second-order flux-integral approach and the discrepancy distribution approach, where we observe in Fig.~\ref{fig:Conservation} a slight accumulation of the round-off error. The third-order flux-integral approach uses the update formulas only for the evolution of the interface values, but the cell average is updated by the explicitly calculated fluxes, which are the same for adjacent cells exactly instead of just up to machine precision. Additionally, we see the correlation between the condition $\sum_{i} \rho_i E_i = 0$ and the conservation of the Momentum and Energy. This condition is directly influenced by the employed Poisson solver and results in different conservation properties if other nodal solvers such as Finite Differences etc. were to be applied.          
	
\section{Conclusion and Outlook}\label{sec:conclusions}

In this paper we have combined the Active Flux method with splitting strategies for the Vlasov--Poisson system. For the original Active Flux method (based on a bubble function) we were able to achieve an overall third-order convergence by combining a higher-order time splitting with a third-order flux integral calculation. The discrepancy distribution formulation, based on pointlike values, allowed a third-order scheme by using only a higher order time splitting method. 

However, two additional observations are important for high-dimensional (6D) Vlasov simulations. First, the Active Flux method is well suited for highly parallel implementations based on domain decomposition due to its very local domain of dependence, at least in its bubble function formulation. The boundary exchange, owing to its high memory requirement for velocity space, is bandwidth-limited and forms a serious bottleneck \cite{schild-raeth-etal:2024}. In this respect, the Active Flux method is even slightly better suited than semi-Lagrangian discontinuous Galerkin methods for massive parallel simulations.

The second and most interesting observation can be seen in Figs.~\ref{fig:convergence_LD}, \ref{fig:convergence_TS}.
In realistic 6D Vlasov simulations, the velocity resolution is limited to a rather small value of $32^3$ up to $48^3$ grid points at most. Figs.~\ref{fig:convergence_LD}, \ref{fig:convergence_TS} show that the discretization error in velocity space dominates over the temporal error introduced by the splitting. The result is not new but confirms the strategy to use high order methods in configuration and velocity space in combination with a relative inexpensive Strang splitting (see e.g. \cite{tanaka-yoshikawa-etal:2017}).

In addition, the long time simulations shown in Figures~\ref{fig:TS_SLD}, \ref{fig:SLD_snapshots}, \ref{fig:TS_snapshots_hist} demonstrate that the Active Flux method is significantly less dissipative than the state-of-the-art PFC method. This is true for all three of the presented splitting approaches. We conclude that the split-step Active Flux method has immense potential for the study of collisionless plasmas. Its low dissipation characteristics and localized formulation make it highly suitable for the development of a massively parallel Vlasov--Maxwell solver to be integrated within the \textit{MuPhy~II} framework~\cite{allmann-rahn-lautenbach-etal:2024,lautenbach-grauer:2018,lautenbach_2024_10547265}, which is a current area of active research in our group. 
	
\section*{Acknowledgements}
We gratefully acknowledge the Gauss Centre for Supercomputing e.V. (www.gauss-centre.eu) for funding this project by providing computing time
through the John von Neumann Institute for Computing (NIC) on the GCS Supercomputer JUWELS \cite{jsc-juwels:2021} at Jülich Supercomputing Centre (JSC).
Computations were conducted on JUWELS/JUWELS-booster and on the DaVinci cluster at TP1 Plasma Research Department. R.G., G.G. and K.K. acknowledge funding from the German Science Foundation DFG through the research unit ``SNuBIC'' (DFG-FOR5409, project ids 463312734 and 530709913).

%% file: main.bbl
\begin{thebibliography}{10}
\expandafter\ifx\csname url\endcsname\relax
  \def\url#1{\texttt{#1}}\fi
\expandafter\ifx\csname urlprefix\endcsname\relax\def\urlprefix{URL }\fi
\expandafter\ifx\csname href\endcsname\relax
  \def\href#1#2{#2} \def\path#1{#1}\fi

\bibitem{palmroth-etal:2018}
M.~Palmroth, U.~Ganse, Y.~Pfau-Kempf, M.~Battarbee, L.~Turc, T.~Brito,
  M.~Grandin, S.~Hoilijoki, A.~Sandroos, S.~von Alfthan, Vlasov methods in
  space physics and astrophysics, Living Reviews in Computational Astrophysics
  4 (2018) 1.

\bibitem{allmann_rahn-lautenbach-grauer:2021}
F.~Allmann-Rahn, S.~Lautenbach, R.~Grauer, {An Energy Conserving Vlasov Solver
  That Tolerates Coarse Velocity Space Resolutions: Simulation of MMS
  Reconnection Events}, Journal of Geophysical Research: Space Physics 127~(2)
  (2022) e2021JA029976.
\newblock \href {https://doi.org/https://doi.org/10.1029/2021JA029976}
  {\path{doi:https://doi.org/10.1029/2021JA029976}}.

\bibitem{nishikawa-etal:2021}
K.~Nishikawa, I.~Du{\c t}an, C.~K{\"o}hn, Y.~Mizuno, {PIC methods in
  astrophysics: simulations of relativistic jets and kinetic physics in
  astrophysical systems}, Living Reviews in Computational Astrophysics 7 (2021)
  1.

\bibitem{leveque:2002}
R.~J. LeVeque, {Finite Volume Methods for Hyperbolic Problems}, Cambridge Texts
  in Applied Mathematics, Cambridge University Press, 2002.

\bibitem{juno-hakim-etal:2018}
J.~Juno, A.~Hakim, J.~TenBarge, E.~Shi, W.~Dorland, Discontinuous {G}alerkin
  algorithms for fully kinetic plasmas, J.\ Comput.\ Phys. 353 (2018) 110--147.

\bibitem{qin-shu:2011}
J.~Qin, C.-W. Shu, Positivity preserving semi-{L}agrangian discontinuous
  {G}alerkin formulation: theoretical analysis and application to the
  {V}lasov-{P}oisson system, J.\ Comput.\ Phys. 230 (2011) 8386--8409.

\bibitem{filbet-sonnendruecker-bertrand:2001}
F.~Filbet, E.~Sonnendr\"ucker, P.~Bertrand, {Conservative Numerical Schemes for
  the {V}lasov Equation}, J.\ Comput.\ Phys. 172 (2001) 166--187.

\bibitem{rossmanith-seal:2011}
J.~Rossmanith, D.~Seal, A positivity-preserving high-order semi-{L}agrangian
  discontinuous {G}alerkin scheme for the {V}lasov-{P}oisson equations, J.\
  Comput.\ Phys. 230 (2011) 6203--6232.

\bibitem{eliasson:2003}
B.~Eliasson, Numerical modelling of the two-dimensional {F}ourier transformed
  {V}lasov-{M}axwell system, J.\ Comput.\ Phys. 190 (2003) 501--522.

\bibitem{delzanno:2015}
G.~Delzanno, Multi-dimensional, fully-implicit, spectral method for the
  {V}lasov-{M}axwell equations with exact conservationlaws in discrete form,
  J.\ Comput.\ Phys. 301 (2015) 338--356.

\bibitem{arber-bennett-etal:2015}
T.~Arber, K.~Bennett, C.~Brady, A.~Lawrence-Douglas, M.~Ramsay, N.~Sircombe,
  P.~Gillies, R.~Evans, H.~Schmitz, A.~Bell, Contemporary particle-in-cell
  approach to laser-plasma modelling, Plasma Physics and Controlled Fusion 57
  (2015) 113001.

\bibitem{masek-gibbon:2010}
M.~Masek, P.~Gibbon, {Mesh-Free Magnetoinduced Plasma Models}, IEEE
  Transactions on Plasma Science 38 (2010) 2377--2382.

\bibitem{cheng-knorr:1976}
C.-Z. Cheng, G.~Knorr, The integration of the {V}lasov equation in
  configuration space, Journal of Computational Physics 22~(3) (1976) 330--351.

\bibitem{schmitz-grauer:2006b}
H.~Schmitz, R.~Grauer, Comparison of time splitting and backsubstitution
  methods for integrating {V}lasov's equation with magnetic fields, Comp.\
  Phys.\ Comm. 175 (2006) 86--92.

\bibitem{roe:2017}
P.~Roe, Multidimensional upwinding, in: R.~Abgrall, C.-W. Shu (Eds.), Handbook
  of Numerical Methods for Hyperbolic Problems, Elsevier, 2017, pp. 53--80.

\bibitem{schmitz-grauer:2006c}
H.~Schmitz, R.~Grauer, Darwin-{V}lasov simulations of magnetised plasmas, J.\
  Comput.\ Phys. 214 (2006) 738--756.

\bibitem{filbet-sonnendruecker:2003}
F.~Filbet, E.~Sonnendr\"ucker, Comparison of {E}ulerian {V}lasov solvers,
  Computer Physics Communications 150 (2003) 247--266.

\bibitem{eymann-roes:2011}
T.~Eymann, P.~Roe, Active {F}lux schemes, in: 49th AIAA Aerospace Sciences
  Meeting including the New Horizons Forum and Aerospace Exposition, 2011, p.
  382.

\bibitem{eymann-roe:2013}
T.~Eymann, P.~Roe, Multidimensional active flux schemes, in: 21st AIAA
  computational fluid dynamics conference (2013).

\bibitem{roe:2021}
P.~Roe, Designing {CFD} methods for bandwidth---a physical approach, Comput. \&
  Fluids 214 (2021) Paper No. 104774, 13.

\bibitem{he2020treatment}
F.~He, P.~L. Roe, {The Treatment of Conservation in the Active Flux method},
  in: AIAA AVIATION 2020 FORUM, 2020, p. 3032.

\bibitem{abgrall-etal:2023b}
R.~Abgrall, W.~Barsukow, {Extensions of Active Flux to arbitrary order of
  accuracy}, ESAIM: M2AN 57~(2) (2023) 991--1027.
\newblock \href {https://doi.org/10.1051/m2an/2023004}
  {\path{doi:10.1051/m2an/2023004}}.

\bibitem{barsukov-holm-etal:2019}
W.~Barsukow, J.~Hohm, C.~Klingenberg, P.~Roe, {The Active Flux Scheme on
  Cartesian Grids and Its Low Mach Number Limit}, J.\ Sci.\ Comput. 81 (2019)
  594--622.

\bibitem{barsukow:2021}
W.~Barsukow, The {A}ctive {F}lux scheme for nonlinear problems, J. Sci. Comput.
  86 (2021) Paper No. 3, 34.

\bibitem{chudzik-helzel-kerkmann:2021}
E.~Chudzik, C.~Helzel, D.~Kerkmann, The {C}artesian grid {A}ctive {F}lux
  method: linear stability and bound preserving limiting, Appl. Math. Comput.
  393 (2021) Paper No.\ 125501, 19.

\bibitem{bai-roe:2021}
Y.~Bai, P.~L. Roe, {Toward Physically-based limiting for the Active Flux
  scheme}, in: AIAA AVIATION 2021 FORUM, 2021, p. 2744.

\bibitem{kiechle2023active}
Y.-F. Kiechle, E.~Chudzik, C.~Helzel, An {A}ctive {F}lux method for the
  {V}lasov-{P}oisson system, in: International Conference on Finite Volumes for
  Complex Applications, Springer, 2023, pp. 93--101.

\bibitem{kiechleNew}
E.~Chudzik, C.~Helzel, Y.~Kiechle, {A Positivity Preserving Active Flux Method
  for the Vlasov--Poisson System}, {personal communication} (2024).

\bibitem{yoshida1990construction}
H.~Yoshida, Construction of higher order symplectic integrators, Physics
  letters A 150~(5-7) (1990) 262--268.

\bibitem{marcowith-ferrand-etal:2020}
A.~Marcowith, G.~Ferrand, M.~Grech, Z.~Meliani, I.~Plotnikov, R.~Walder,
  Multi-scale simulations of particle acceleration in astrophysical systems,
  Living Reviews in Computational Astrophysics 6~(1) (2020) 1.
\newblock \href {https://doi.org/10.1007/s41115-020-0007-6}
  {\path{doi:10.1007/s41115-020-0007-6}}.

\bibitem{kormann-reuter-rampp:2019}
K.~Kormann, K.~Reuter, M.~Rampp, A massively parallel semi-{L}agrangian solver
  for the six-dimensional {V}lasov–{P}oisson equation, The International
  Journal of High Performance Computing Applications 33~(5) (2019) 924--947.
\newblock \href {https://doi.org/10.1177/1094342019834644}
  {\path{doi:10.1177/1094342019834644}}.

\bibitem{chudzik-etal:2024}
E.~Chudzik, C.~Helzel, M.~Lukáčová-Medvid’ová, Active {Flux} {Methods}
  for {Hyperbolic} {Systems} {Using} the {Method} of {Bicharacteristics},
  Journal of Scientific Computing 99~(1) (2024) 16.
\newblock \href {https://doi.org/10.1007/s10915-024-02462-z}
  {\path{doi:10.1007/s10915-024-02462-z}}.

\bibitem{abgrall-etal:2023}
R.~Abgrall, W.~Barsukow, C.~Klingenberg,
  \href{https://arxiv.org/abs/2310.00683}{{The Active Flux method for the Euler
  equations on Cartesian grids}} (2023).
\newblock \href {http://arxiv.org/abs/2310.00683} {\path{arXiv:2310.00683}}.
\newline\urlprefix\url{https://arxiv.org/abs/2310.00683}

\bibitem{maeng2017advective}
J.~Maeng, On the advective component of {A}ctive {F}lux schemes for nonlinear
  hyperbolic conservation laws, Ph.D. thesis (2017).

\bibitem{robidoux2008polynomial}
N.~Robidoux, Polynomial histopolation, superconvergent degrees of freedom, and
  pseudospectral discrete hodge operators, Unpublished: http://www. cs.
  laurentian. ca/nrobidoux/prints/super/histogram. pdf (2008).

\bibitem{kormann2024dual}
K.~Kormann, E.~Sonnendr{\"u}cker, A dual grid geometric electromagnetic
  particle in cell method, SIAM Journal on Scientific Computing 46~(5) (2024)
  B621--B646.

\bibitem{schild-raeth-etal:2024}
N.~Schild, M.~Räth, S.~Eibl, K.~Hallatschek, K.~Kormann, A performance
  portable implementation of the semi-{L}agrangian algorithm in six dimensions,
  Computer Physics Communications 295 (2024) 108973.
\newblock \href {https://doi.org/10.1016/j.cpc.2023.108973}
  {\path{doi:10.1016/j.cpc.2023.108973}}.

\bibitem{sonnendrucker2013numerical}
E.~Sonnendr{\"u}cker, K.~Kormann, Numerical methods for vlasov equations,
  Lecture notes 107 (2013) 108.

\bibitem{tanaka-yoshikawa-etal:2017}
S.~Tanaka, K.~Yoshikawa, T.~Minoshima, N.~Yoshida, {Multidimensional
  Vlasov–Poisson Simulations with High-order Monotonicity- and
  Positivity-preserving Schemes}, The Astrophysical Journal 849~(2) (2017) 76.
\newblock \href {https://doi.org/10.3847/1538-4357/aa901f}
  {\path{doi:10.3847/1538-4357/aa901f}}.

\bibitem{allmann-rahn-lautenbach-etal:2024}
F.~Allmann-Rahn, S.~Lautenbach, M.~Deisenhofer, R.~Grauer, The muphyii code:
  Multiphysics plasma simulation on large hpc systems, Computer Physics
  Communications 296 (2024) 109064.
\newblock \href {https://doi.org/10.1016/j.cpc.2023.109064}
  {\path{doi:10.1016/j.cpc.2023.109064}}.

\bibitem{lautenbach-grauer:2018}
S.~Lautenbach, R.~Grauer, Multiphysics simulations of collisionless plasmas,
  Frontiers in Physics 6 (2018).
\newblock \href {https://doi.org/10.3389/fphy.2018.00113}
  {\path{doi:10.3389/fphy.2018.00113}}.

\bibitem{lautenbach_2024_10547265}
S.~Lautenbach, F.~Allmann-Rahn, R.~Grauer, muphy2 (Jan. 2024).
\newblock \href {https://doi.org/10.5281/zenodo.10547265}
  {\path{doi:10.5281/zenodo.10547265}}.

\bibitem{jsc-juwels:2021}
D.~Alvarez, {JUWELS Cluster and Booster: Exascale Pathfinder with Modular
  Supercomputing Architecture at Juelich Supercomputing Centre} 7  A183--A183.
\newblock \href {https://doi.org/10.17815/jlsrf-7-183}
  {\path{doi:10.17815/jlsrf-7-183}}.

\end{thebibliography}
